\documentclass[10pt,a4paper]{article}
\usepackage[margin=2cm]{geometry}
\usepackage{amsmath,amsthm,amsfonts,amssymb,amscd,cite,mathrsfs,graphicx,url}
\usepackage{latexsym}
\usepackage{enumitem}
\usepackage{thmtools}
\usepackage[usenames, dvipsnames]{color}
\definecolor{mygray}{gray}{0.7}

\usepackage[table]{xcolor}
\usepackage{wasysym,upgreek}
\usepackage{mathdots}
\usepackage[enableskew]{youngtab}
\usepackage{tikz}
\usetikzlibrary{tikzmark,decorations.pathreplacing}
\usetikzlibrary{arrows,matrix}
\usetikzlibrary{positioning}
\usetikzlibrary{decorations}
\usepackage[all,cmtip]{xy}
\usepackage{hyperref}
\usepackage{comment}
\usepackage{diagbox}
\usepackage{titlesec}
\usepackage{bbold}
\titleformat{\sec}
{\normalfont\fontsize{16}{20}\bfseries}{\thesec}{1em.}{}
\definecolor{amber}{rgb}{1.0,0.75,0.0}
\definecolor{applegreen}{rgb}{0.55,0.71,0.0}
\definecolor{byzantium}{rgb}{0.44, 0.16, 0.39}
\definecolor{cadmiumorange}{rgb}{0.93, 0.53, 0.18}
\definecolor{darkcyan}{rgb}{0.0, 0.55, 0.55}

\newtheorem{conjecture}{Conjecture}[section]
\newtheorem{corollary}{Corollary}[section]
\newtheorem{lemma}{Lemma}[section]
\newtheorem{definition}{Definition}[section]
\newtheorem{theorem}{Theorem}[section]

\newtheorem{example}{Example}[section]

\newtheorem{rmk}{Remark}[section]

\allowdisplaybreaks[4]

\let\oldbibliography\thebibliography
\renewcommand{\thebibliography}[1]{%
  \oldbibliography{#1}%
  \setlength{\itemsep}{-2pt}%
}
\usepackage{fancyhdr}
\fancypagestyle{mypagestyle}{%
  \fancyhf{}
  \fancyhead[OC]{Jean-Luc Baril$^\dag$, Nathana\"el Hassler$^\dag$, Sergey Kirgizov$^\dag$ and José L. Ram\ii rez$^\ddag$}
  \fancyfoot[C]{{\it ECA}  {\bf x:x} (xxxx) Article \#SxRx\hfill\arabic{page}}%
}
\newcommand{\seqnum}[1]{\href{http://oeis.org/#1}{\underline{#1}}}

\makeatletter
\def\BState{\State\hskip-\ALG@thistlm}
\makeatother
\def\a{0.5cm}

\usepackage[colorinlistoftodos]{todonotes}
 
 \newcommand\scalemath[2]{\scalebox{#1}{\mbox{\ensuremath{\displaystyle #2}}}}

\def\N{\mathbb{N}}
\def\Z{\mathbb{Z}}

\def\Eu{(4,2)}
\def\Nu{(2,4)}
\def\Ed{(4,-2)}
\def\Nd{(2,-4)}
 \def\ii{\'i}

\newcommand{\Mod}[1]{\ (\mathrm{mod}\ #1)}

\begin{document}
\baselineskip=0.20in

\makebox[\textwidth]{%
\hglue-15pt
\begin{minipage}{0.6cm}	
\vskip9pt
\end{minipage} \vspace{-\parskip}
\begin{minipage}[t]{11cm}
\footnotesize{ {\bf Enumerative Combinatorics and Applications} \\ \underline{ecajournal.haifa.ac.il}}
\end{minipage}
\hfill
\begin{minipage}[t]{5.4cm}
\normalsize {\it ECA}  {\bf X:X} (xxxx) Article \#SxRx
\end{minipage}}
\vskip36pt

\begin{center}
{\large \bf Grand zigzag Knight's paths}\\[10pt]

Jean-Luc Baril$^\dag$, Nathana\"el Hassler$^\dag$, Sergey Kirgizov$^\dag$ and José L. Ram\ii rez$^\ddag$\\[20pt]

\footnotesize {\it $^\dag$ LIB, Universit\'e de Bourgogne,   B.P. 47 870, 21078, Dijon Cedex, France\\
Email: barjl@u-bourgogne.fr, nathanael.hassler@ens-rennes.fr, sergey.kirgizov@u-bourgogne.fr}\\[10pt]

\footnotesize {\it $^\ddag$ Departamento de Matem\'aticas,  Universidad Nacional de Colombia,  Bogot\'a, Colombia\\
Email: jlramirezr@unal.edu.co}\\[20pt]

\noindent {\footnotesize {\bf Received}: April XX, xxxx},
{\footnotesize {\bf Accepted}: August XX, xxxx}, {\footnotesize {\bf Published}: August XX, xxxx}\\

\noindent The authors: Released under the CC BY-ND license (International 4.0)
\end{center}

\setcounter{page}{1} \thispagestyle{empty}

\baselineskip=0.30in

\normalsize

\noindent
{\sc Abstract:} We study the enumeration of different classes of grand knight's paths in the plane.  In particular, we focus on the subsets of zigzag knight's paths that are subject to constraints.  These constraints include ending at $y$-coordinate  0,  bounded by a horizontal line, confined within a tube, among other considerations. 
We present our results using  generating functions or direct closed-form expressions. We derive asymptotic results, finding approximations for quantities such as the probability that a zigzag knight's path stays in some area of the plane, or for the average of the altitude of such a path. Additionally, we exhibit some bijections between grand zigzag knight's paths and some pairs of compositions.
\bigskip



\section{Introduction}

In combinatorics, the enumeration of different classes of lattice paths according to several parameters is often studied. One of the motivations is to exhibit new one-to-one correspondences between these classes and various objects from  other domains, such as computer science, biology, and physics~\cite{Sta}. For instance, they have close connections with RNA structures, pattern-avoiding permutations, directed animals, and other related topics~\cite{Bar1,Knu,Sta}. We refer to \cite{Ban,Barc,Bar,Deu,Flo, Man1,Mer,Sap,Sun,SunZha}, and the references therein, for such studies. 

Labelle and Yeh \cite{Lab}  investigate knight's paths, that are lattice paths in $\mathbb{N}^2$ that start at the origin and consist of steps $N=(1,2)$, $\bar{N}=(1,-2)$, $E=(2,1)$, and $\bar{E}=(2,-1)$. Note that these steps correspond to the right moves of a knight on a chessboard. Recently, in \cite{BarRa}, the authors focus on zigzag knight's paths, i.e.,  knight's paths where the direction of the steps alternate up and down, or equivalently knight's paths avoiding  the consecutive patterns $NN$, $NE$, $EN$, $EE$, $\bar{N}\bar{N}$, $\bar{N}\bar{E}$, $\bar{E}\bar{N}$, $\bar{E}\bar{E}$. They prove that such paths ending on the $x$-axis with a given number of steps are enumerated by the well-known Catalan numbers (see \href{https://oeis.org/A000108}{A000108} in Sloane's On-line Encyclopedia  of  Integer Sequences~\cite{OEIS}), and they exhibit a bijection between these paths and Dyck paths (lattice paths in $\mathbb{N}^2$ starting at the origin, ending on the $x$-axis and made of steps $U=(1,1)$ and $D=(1,-1)$).

In this work, we extend the aforementioned study by allowing knight's paths to go below the $x$-axis. 
\begin{definition} A {\it grand knight's path}  is a lattice path in $\mathbb{Z}^2$ that starts at the origin and consists of steps $N=(1,2)$, $\bar{N}=(1,-2)$, $E=(2,1)$, and $\bar{E}=(2,-1)$.    
\end{definition}
 The \textit{size} of such a path is the $x$-coordinate of its last point. The \emph{empty path} $\epsilon$ is a  path of size~$0$.  The {\it height} is the maximal $y$-coordinate reached by a point of the path, and the {\it altitude} is the $y$-coordinate of the last point of the path. We also define grand zigzag knight's paths as follows.
\begin{definition}\label{d12} A {\it grand zigzag knight’s path} is a grand knight’s path with the
additional property that the vertical components of two consecutive steps cannot be in the same direction, i.e., two
consecutive steps cannot be $NN$, $NE$, $\bar{N}\bar{N}$, $\bar{N}\bar{E}$, $EE$, $EN$, $\bar{E}\bar{E}$, $\bar{E}\bar{N}$.   
\end{definition}

\begin{figure}[ht]
 \begin{center}
        \begin{tikzpicture}[scale=0.15]
    \draw (\a,\a)-- ++(38,0);   \filldraw (\a,\a) circle (8pt);
    \draw[dashed,line width=0.1mm] (\a,2.5)-- ++(38,0);
    \draw[dashed,line width=0.1mm] (\a,4.5)-- ++(38,0);
    \draw[dashed,line width=0.1mm] (\a,6.5)-- ++(38,0);
    \draw[dashed,line width=0.1mm] (\a,-1.5)-- ++(38,0);
    \draw[dashed,line width=0.1mm] (\a,-3.5)-- ++(38,0);
    \draw[dashed,line width=0.1mm] (\a,-5.5)-- ++(38,0);
    \draw (\a,-5.5) -- (\a,6.5);
    \draw[dashed,line width=0.1mm] (2.5,-5.5) --(2.5,6.5);
    \draw[dashed,line width=0.1mm] (4.5,-5.5) --(4.5,6.5);
    \draw[dashed,line width=0.1mm] (6.5,-5.5) --(6.5,6.5);
    \draw[dashed,line width=0.1mm] (8.5,-5.5) --(8.5,6.5);
    \draw[dashed,line width=0.1mm] (10.5,-5.5) --(10.5,6.5);
    \draw[dashed,line width=0.1mm] (12.5,-5.5) --(12.5,6.5);
    \draw[dashed,line width=0.1mm] (14.5,-5.5) --(14.5,6.5);
    \draw[dashed,line width=0.1mm] (16.5,-5.5) --(16.5,6.5);
    \draw[dashed,line width=0.1mm] (18.5,-5.5) --(18.5,6.5);
    \draw[dashed,line width=0.1mm] (20.5,-5.5) --(20.5,6.5);
    \draw[dashed,line width=0.1mm] (22.5,-5.5) --(22.5,6.5);
    \draw[dashed,line width=0.1mm] (24.5,-5.5) --(24.5,6.5);
    \draw[dashed,line width=0.1mm] (26.5,-5.5) --(26.5,6.5);
    \draw[dashed,line width=0.1mm] (28.5,-5.5) --(28.5,6.5);
    \draw[dashed,line width=0.1mm] (30.5,-5.5) --(30.5,6.5);
    \draw[dashed,line width=0.1mm] (32.5,-5.5) --(32.5,6.5);
    \draw[dashed,line width=0.1mm] (34.5,-5.5) --(34.5,6.5);
    \draw[dashed,line width=0.1mm] (36.5,-5.5) --(36.5,6.5);
    \draw[dashed,line width=0.1mm] (38.5,-5.5) --(38.5,6.5);
    \draw[solid,line width=0.4mm] (\a,\a)-- ++\Eu  -- ++\Nu-- ++\Ed-- ++\Ed--++\Nd-- ++\Ed-- ++\Eu-- ++\Ed-- ++\Nu-- ++\Eu-- ++\Nu-- ++\Nd;
    \filldraw (\a,\a)++\Eu circle (8pt);\filldraw (\a,\a)++\Eu++\Nu circle (8pt);
    \filldraw (\a,\a)++\Eu++\Nu++\Ed circle (8pt);
    \filldraw (\a,\a)++\Eu++\Nu++\Ed++\Ed circle (8pt);
    \filldraw (\a,\a)++\Eu++\Nu++\Ed++\Ed++\Nd circle (8pt);
     \filldraw (\a,\a)++\Eu++\Nu++\Ed++\Ed++\Nd++\Ed circle (8pt);
     \filldraw (\a,\a)++\Eu++\Nu++\Ed++\Ed++\Nd++\Ed++\Eu circle (8pt);
     \filldraw (\a,\a)++\Eu++\Nu++\Ed++\Ed++\Nd++\Ed++\Eu++\Ed circle (8pt);
     \filldraw (\a,\a)++\Eu++\Nu++\Ed++\Ed++\Nd++\Ed++\Eu++\Ed++\Nu circle (8pt);
      \filldraw (\a,\a)++\Eu++\Nu++\Ed++\Ed++\Nd++\Ed++\Eu++\Ed++\Nu++\Eu circle (8pt);
      \filldraw (\a,\a)++\Eu++\Nu++\Ed++\Ed++\Nd++\Ed++\Eu++\Ed++\Nu++\Eu++\Nu circle (8pt);
       \filldraw (\a,\a)++\Eu++\Nu++\Ed++\Ed++\Nd++\Ed++\Eu++\Ed++\Nu++\Eu++\Nu++\Nd circle (8pt);
\end{tikzpicture}\quad 
 \begin{tikzpicture}[scale=0.15]
    \draw (\a,\a)-- ++(38,0);
    \draw[dashed,line width=0.1mm] (\a,2.5)-- ++(38,0);
    \draw[dashed,line width=0.1mm] (\a,4.5)-- ++(38,0);
    \draw[dashed,line width=0.1mm] (\a,6.5)-- ++(38,0);
    \draw[dashed,line width=0.1mm] (\a,-1.5)-- ++(38,0);
    \draw[dashed,line width=0.1mm] (\a,-3.5)-- ++(38,0);
    \draw[dashed,line width=0.1mm] (\a,-5.5)-- ++(38,0);
    \draw (\a,-5.5) -- (\a,6.5);
    \draw[dashed,line width=0.1mm] (2.5,-5.5) --(2.5,6.5);
    \draw[dashed,line width=0.1mm] (4.5,-5.5) --(4.5,6.5);
    \draw[dashed,line width=0.1mm] (6.5,-5.5) --(6.5,6.5);
    \draw[dashed,line width=0.1mm] (8.5,-5.5) --(8.5,6.5);
    \draw[dashed,line width=0.1mm] (10.5,-5.5) --(10.5,6.5);
    \draw[dashed,line width=0.1mm] (12.5,-5.5) --(12.5,6.5);
    \draw[dashed,line width=0.1mm] (14.5,-5.5) --(14.5,6.5);
    \draw[dashed,line width=0.1mm] (16.5,-5.5) --(16.5,6.5);
    \draw[dashed,line width=0.1mm] (18.5,-5.5) --(18.5,6.5);
    \draw[dashed,line width=0.1mm] (20.5,-5.5) --(20.5,6.5);
    \draw[dashed,line width=0.1mm] (22.5,-5.5) --(22.5,6.5);
    \draw[dashed,line width=0.1mm] (24.5,-5.5) --(24.5,6.5);
    \draw[dashed,line width=0.1mm] (26.5,-5.5) --(26.5,6.5);
    \draw[dashed,line width=0.1mm] (28.5,-5.5) --(28.5,6.5);
    \draw[dashed,line width=0.1mm] (30.5,-5.5) --(30.5,6.5);
    \draw[dashed,line width=0.1mm] (32.5,-5.5) --(32.5,6.5);
    \draw[dashed,line width=0.1mm] (34.5,-5.5) --(34.5,6.5);
    \draw[dashed,line width=0.1mm] (36.5,-5.5) --(36.5,6.5);
    \draw[dashed,line width=0.1mm] (38.5,-5.5) --(38.5,6.5);
    \draw[solid,line width=0.4mm] (\a,\a)-- ++\Eu  -- ++\Nd-- ++\Eu-- ++\Ed--++\Eu-- ++\Nd-- ++\Eu-- ++\Ed-- ++\Eu-- ++\Nd-- ++\Eu;
    \filldraw (\a,\a)++\Eu circle (8pt); \filldraw (\a,\a) circle (8pt);
    \filldraw (\a,\a)++\Eu++\Nd circle (8pt);
    \filldraw (\a,\a)++\Eu++\Nd++\Eu circle (8pt);
    \filldraw (\a,\a)++\Eu++\Nd++\Eu++\Ed circle (8pt);
    \filldraw (\a,\a)++\Eu++\Nd++\Eu++\Ed++\Eu circle (8pt);
     \filldraw (\a,\a)++\Eu++\Nd++\Eu++\Ed++\Eu++\Nd circle (8pt);
     \filldraw (\a,\a)++\Eu++\Nd++\Eu++\Ed++\Eu++\Nd++\Eu circle (8pt);
     \filldraw (\a,\a)++\Eu++\Nd++\Eu++\Ed++\Eu++\Nd++\Eu++\Ed circle (8pt);
     \filldraw (\a,\a)++\Eu++\Nd++\Eu++\Ed++\Eu++\Nd++\Eu++\Ed++\Eu circle (8pt);
      \filldraw (\a,\a)++\Eu++\Nd++\Eu++\Ed++\Eu++\Nd++\Eu++\Ed++\Eu++\Nd circle (8pt);
      \filldraw (\a,\a)++\Eu++\Nd++\Eu++\Ed++\Eu++\Nd++\Eu++\Ed++\Eu++\Nd++\Eu circle (8pt);
       \filldraw (\a,\a)++\Eu++\Nd++\Eu++\Ed++\Eu++\Nd++\Eu++\Ed++\Eu++\Nd++\Eu circle (8pt);
\end{tikzpicture}
               \end{center}
         \caption{On the left: a grand knight's path of size 19, height 3, and altitude~1. On the right: a grand zigzag knight's path of size 19, height 1, and altitude~$-2$.}
         \label{figg1}
\end{figure}
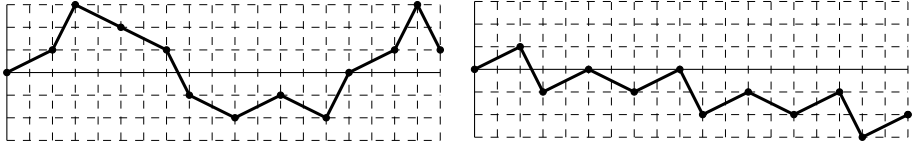

See Figure~\ref{figg1} for an illustration of two grand knight's paths,  with the second one also being a grand zigzag knight's path.

The enumeration of grand knight's paths of a given size can be obtained using the work of Banderier and Flajolet~\cite{bf}. 
Indeed, it suffices to consider the step polynomial $P(z,u)=zu^2+\frac{z}{u^2}+\frac{z^2}{u}+z^2u$, where $z$ (resp.~$u$) marks the size (resp.~height) of the step, i.e., $zu^2$, $\frac{z}{u^2}$, $\frac{z^2}{u}$, and $z^2u$ correspond to the steps $N$, $\bar{N}$, $\bar{E}$, and $E$, respectively. The Laurent series providing the number of grand knight's paths with respect to the size (marked with $z$), the number of steps (marked with $x$), and the altitude (marked with $u$) is given by 
$$\frac{1}{1-xP(z,u)}.$$ 
Fixing $x=u=1$, we deduce the generating function $1/(1-2z-2z^2)$ for the number of grand knight's paths with respect to the size. The coefficient of $z^n$ is  $\sum_{k=0}^{n+1} \binom{n+1}{2k + 1}3^k$, which corresponds to the sequence   \seqnum{A002605}  in the OEIS \cite{OEIS}. The first values of $u_n$ are 
$$1,\ 2,\ 6,\ 16,\ 44,\ 120,\ 328,\ 896,\ 2448,\ 6688,\ 18272.$$ An asymptotic approximation  is given by $\sqrt{3}/6\cdot (\sqrt{3}+1)^{n+1}$.

The generating function for the number of grand knight's paths ending on the $x$-axis 
is 
$$x\left(\frac{\mathbf{u}_3(x)'}{\mathbf{u}_3(x)}
+\frac{\mathbf{u}_4(x)'}{\mathbf{u}_4(x)}\right)\bigg|_{x=1},$$ where $\mathbf{u}_3(x)$ and $\mathbf{u}_4(x)$ are the small roots in $u$ of $1-xP(u)$ (i.e. the roots $\mathbf{u}_i(x)$ that satisfy $\lim_{x \to 0} \mathbf{u}_i(x) = 0$), and the derivative with respect to $x$.
The bivariate generating function for the number of grand knight's paths of a given size and altitude $k\geq 0$ is then 
$$W_k(z,u)=\left(  \frac{x}{k}  \left(
\mathbf{u}_3(x)^{-k}\right)' +
\frac{x}{k}\left(\mathbf{u}_4(x)^{-k}\right)'
\right)
\bigg|_{x=1}.$$

In the subsequent sections of the paper, among other results, we recover the aforementioned formulae, although in a different statement (without involving any derivatives).

Recently, 
Asinowski, Bacher, Banderier, Gittenberger, and Roitner~\cite{abbcg1, abbcg2, abbcg3, asinowski-banderier-roitner, ab} developed 
a method for solving enumerative problems that arise out of paths avoiding consecutive patterns. 
Grand zigzag knight’s paths could be modeled as grand knight’s paths avoiding eight consecutive patterns of size 2. 
However, in this paper we present a more specific method, using the point of view of Prodinger \cite{Pro} based on the kernel method. It consists in studying two generating functions that count the desired paths by distinguishing the direction of the last step (up or down). Then, we provide formulae for different kinds of constraints, give bijections, direct closed-form expressions for the $n$-th coefficient, and asymptotic results.

{\bf Outline of the paper.} In Section~2, we apply the kernel method for providing generating function of the number of grand knight's paths with respect to the size and the altitude. We deduce asymptotic approximations for the number of grand knight's paths of a given size and a non-negative altitude, and for the expected altitude of a grand knight's path ending at a positive altitude. In the following sections, we focus on grand zigzag knight's paths. In Section~3, we give the counterpart of Section~2 for grand zigzag knight's paths. Moreover, we exhibit a bijection between pairs of integer compositions and grand zigzag knight's paths ending at $(n,k)$, where $n$ and $k$ have the same parity. As a byproduct, we obtain a closed form for the number of grand zigzag knight's paths ending at $(n,k)$, and the expected value for the number of steps of a grand zigzag knight's path ending at $(n,k)$; asymptotic approximations are also derived. In Section~4, we consider grand zigzag knight's paths staying in a region delimited by horizontal lines. We make a similar study as for the previous sections, which allows us to estimate  asymptotic approximations for the probability that a grand zigzag knight's path chosen uniformly at random among all grand zigzag knight's paths stays above a horizontal line. We end by giving an appendix for the general case of grand zigzag knight's paths staying in a general tube.

\section{Unrestricted grand knight's paths}

In this section, we enumerate grand knight's paths starting at the origin $(0,0)$ and ending at $(n,k)$ for $n\geq 0$, $k\in\mathbb{Z}$. Let $\mathcal{H}_{n,k}$ be the set of such paths and let $h_{n,k}$ be its cardinality.  Due to the symmetry with respect to the $x$-axis, we obviously have $h_{n,k}=h_{n,-k}$. Therefore it is sufficient to focus on non-negative $k$. For (fixed) $k\geq 0$, let $h_k=\sum_{n\geq 0}h_{n,k} z^n$ be the generating function for the number of grand knight's paths of altitude $k$ with respect to the size. A non-empty grand knight's path of altitude $k$ and size $n\geq 1$ is of one of the following forms ($a$) $P\bar{N}$ with $P\in\mathcal{H}_{n-1,k+2}$, ($b$) $P\bar{E}$ with $P\in\mathcal{H}_{n-2,k+1}$, ($c$) $PE$ with $P\in\mathcal{H}_{n-2,k-1}$, ($d$) $PN$ with $P\in\mathcal{H}_{n-1,k-2}$. Therefore we easily obtain the following equations after considering the symmetry $h_{-k}=h_k$: 

\begin{align}
    h_0&=1+2zh_2+2z^2h_1,\label{h_0}\\
     h_1&=z^2(h_0+h_2)+z(h_1+h_3), \label{h_1}\\
    h_k&=z^2(h_{k-1}+h_{k+1})+z(h_{k-2}+h_{k+2}),  \text{ for } k\geq 2. \label{h_k}
\end{align}
Let $H(u,z)=\sum_{k\geq 0}h_k(z)u^k$ be the bivariate generating function for the number of grand knight's paths with respect to the size and the altitude. We will write $H(u)$ for short. Thus, by multiplying \eqref{h_k} by $u^{k}$, summing over $k\geq2$, with $h_0+h_1u$ we obtain
\begin{align*}
    H(u)=h_0+&h_1u+z^2\sum_{k\geq 2} (h_{k-1}+h_{k+1})u^k +z\sum_{k\geq 2} (h_{k-2}+h_{k+2})u^k,
\end{align*}
and 
\begin{equation}\label{equaH}
  \begin{split}
     H(u) =h_0+h_1u+z^2\left(u(H(u)-h_0)+\frac{H(u)-h_0-h_1u-h_2u^2}{u}\right)\\
    +z\left( u^2H(u)+\frac{H(u)-h_0-h_1u-h_2u^2-h_3u^3}{u^2}\right). 
  \end{split} 
\end{equation}
From (\ref{h_0}) and (\ref{h_1}) we obtain $h_2=-\frac{1}{2z}+\frac{h_0}{2z}-zh_1$ and $h_3=\frac{1}{2}+(z^2+\frac{1}{z}-1)h_1-(\frac{1}{2}+z)h_0$. By replacing $h_2$ and $h_3$ by these expressions in (\ref{equaH}), and after simplifying, we eventually have

\begin{equation}\label{main u^2}
    H(u)\left(u^2-zu^4-z-z^2u-z^2u^3\right)=\frac{u^2}{2}+h_0\left(\frac{u^2}{2}-z-z^2u\right)+zuh_1(u^2-1).
\end{equation}

 We use the kernel method (see \cite{Pro}) to determine $h_0$ and $h_1$. Let $K(u)=u^2-zu^4-z-z^2u-z^2u^3$ be the \textit{kernel} of Equation~(\ref{main u^2}),
 and let $u_1$ and $u_2$ be two of the roots of  $K(u)$:
\begin{equation}\label{equ1}
u_1=\frac{-z^2+\sqrt{z^4+8z^2+4z}}{4z}+\frac{1}{2\sqrt{2}}\sqrt{\frac{z^3-z\sqrt{z^4+8z^2+4z}-4z+2}{z}}
\end{equation}
and
\begin{equation}\label{equ2}
u_2=\frac{-z^2-\sqrt{z^4+8z^2+4z}}{4z}-\frac{1}{2\sqrt{2}}\sqrt{\frac{z^3+z\sqrt{z^4+8z^2+4z}-4z+2}{z}}.
\end{equation}
Note that $u_1=\mathbf{u}_1(1)$ and $u_2=\mathbf{u}_2(1)$
where $\mathbf{u}_1(x)$ and $\mathbf{u}_2(x)$
are the two other roots of the step polynomial $1-xP(z,u)$ (see Introduction), which implies that $1/(u-u_1)$ and $1/(u-u_2)$ have no power series expansions around $(u,z)=(0,0)$. This implies that $(u-u_1)(u-u_2)$ must be a factor of the numerator. Then by (\ref{main u^2}), we have  
$$\frac{u_i^2}{2}+h_0\left(\frac{u_i^2}{2}-z-z^2u_i\right)+zu_ih_1(u_i^2-1)=0,  \quad \text{for } i=1,2.$$
After solving for $h_0$ and $h_1$, we obtain

\begin{equation}\label{eqh0}
h_0=\frac{1+u_1u_2}{1+u_1u_2-2z^2(u_1+u_2)-2z+2z(1-u_1^2-u_2^2)u_1^{-1}u_2^{-1}}
\end{equation}
and 
\begin{equation}\label{eqh1}
h_1=\frac{u_1+u_2}{1+u_1u_2}h_0.
\end{equation}

Then we can deduce $H(u)$ from (\ref{main u^2}), which yields the following.

\begin{theorem}
The bivariate generating function for grand knight's paths with respect to  the size and the altitude is
    $$H(u)=\frac{-h_1 u +h_0u_1 u_2}{\left(u -u_1 \right) \left(u -u_2 \right)},$$
    where $u_1, u_2, h_0,h_1$ are defined in (\ref{equ1}),~(\ref{equ2}),~(\ref{eqh0}) and (\ref{eqh1}). 
\end{theorem}

By decomposing $H(u)$ into partial fractions, we deduce a close form for the coefficient $[u^k]H(u)$ of $u^k$ in $H(u)$. 

\begin{corollary}
    The generating function for grand knight's paths of altitude $k$ with respect to the size is
    $$[u^k]H(u)=\frac{h_0u_1-h_1}{u_1-u_2}u_2^{-k}- \frac{h_0u_2-h_1}{u_1-u_2}u_1^{-k}.$$
\end{corollary}

    The following table presents the numbers $h_{n,k}$ for $0\leq n\leq 9$, $0\leq k\leq 9$. It is given in the entry \href{https://oeis.org/A096608}{A096608}  of \cite{OEIS}.
    
\footnotesize$$\left(\begin{array}{ccccccccccc}
    1 & 0 & 0 & 0 & 0 & 0 & 0 & 0 & 0 & 0 &\cdots\\
    0 & 0 & 1 & 0 & 0 & 0 & 0 & 0 & 0 & 0&\cdots\\
    2 & 1 & 0 & 0 & 1 & 0 & 0 & 0 & 0 & 0&\cdots\\
    0 & 2 & 3 & 2 & 0 & 0 & 1 & 0 & 0 & 0&\cdots\\
    \mathbf{8} & 6 & 1 & 3 & 4 & 3 & 0 & 0 & 1 & 0&\cdots\\
    6 & 12 & 16 & 12 & 3 & 4 & 5 & 4 & 0 & 0&\cdots\\
    44 & 33 & 18 & 21 & 27 & 20 & 6 & 5 & 6 & 5&\cdots\\
    60 & 76 & 95 & 72 & 40 & 34 & 41 & 30 & 10 & 6&\cdots\\
    256 & 210 & 154 & 155 & 177 & 135 & 75 & 52 & 58 & 42&\cdots\\
    460 & 520 & 581 & 480 & 335 & 288 & 299 & 228 & 126 & 76&\cdots\\
 \vdots  &  \vdots& \vdots & \vdots & \vdots & \vdots& \vdots  &  \vdots & \vdots & \vdots&\ddots
\end{array}\right).$$
\normalsize

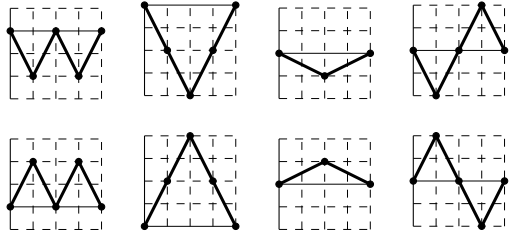
\begin{figure}[ht]
  \begin{center}
\begin{tikzpicture}[scale=0.15]
    \draw[line width=0.1mm] (\a,\a)-- ++(8,0);
    \draw[dashed,line width=0.1mm] (\a,2.5)-- ++(8,0);
    \draw[dashed,line width=0.1mm] (\a,-1.5)-- ++(8,0);
    \draw[dashed,line width=0.1mm] (\a,-3.5)-- ++(8,0);
    \draw[dashed,line width=0.1mm] (\a,-5.5)-- ++(8,0);
    \draw (\a,-5.5) -- (\a,2.5);
    \draw[dashed,line width=0.1mm] (2.5,-5.5) --(2.5,2.5);
    \draw[dashed,line width=0.1mm] (4.5,-5.5) --(4.5,2.5);
    \draw[dashed,line width=0.1mm] (6.5,-5.5) --(6.5,2.5);
    \draw[dashed,line width=0.1mm] (8.5,-5.5) --(8.5,2.5);
    \draw[solid,line width=0.4mm] (\a,\a)-- ++\Nd -- ++\Nu -- ++\Nd -- ++\Nu;
    \filldraw (\a,\a) circle (8pt);
    \filldraw (\a,\a)++\Nd circle (8pt);
    \filldraw (\a,\a)++\Nd++\Nu circle (8pt);
    \filldraw (\a,\a)++\Nd++\Nu++\Nd circle (8pt);
    \filldraw (\a,\a)++\Nd++\Nu++\Nd++\Nu circle (8pt);
    \end{tikzpicture}
\quad
\begin{tikzpicture}[scale=0.15]
    \draw[line width=0.1mm] (\a,\a)-- ++(8,0);
    \draw[dashed,line width=0.1mm] (\a,-1.5)-- ++(8,0);
    \draw[dashed,line width=0.1mm] (\a,-3.5)-- ++(8,0);
    \draw[dashed,line width=0.1mm] (\a,-5.5)-- ++(8,0);
    \draw[dashed,line width=0.1mm] (\a,-7.5)-- ++(8,0);
    \draw (\a,-7.5) -- (\a,\a);
    \draw[dashed,line width=0.1mm] (2.5,-7.5) --(2.5,\a);
    \draw[dashed,line width=0.1mm] (4.5,-7.5) --(4.5,\a);
    \draw[dashed,line width=0.1mm] (6.5,-7.5) --(6.5,\a);
    \draw[dashed,line width=0.1mm] (8.5,-7.5) --(8.5,\a);
    \draw[solid,line width=0.4mm] (\a,\a)-- ++\Nd -- ++\Nd -- ++ \Nu -- ++ \Nu;
     \filldraw (\a,\a) circle (8pt);
    \filldraw (\a,\a)++\Nd circle (8pt);
    \filldraw (\a,\a)++\Nd++\Nd circle (8pt);
    \filldraw (\a,\a)++\Nd++\Nd++\Nu circle (8pt);
    \filldraw (\a,\a)++\Nd++\Nd++\Nu++\Nu circle (8pt);
\end{tikzpicture}
\quad
\begin{tikzpicture}[scale=0.15]
    \draw[line width=0.1mm] (\a,\a)-- ++(8,0);
    \draw[dashed,line width=0.1mm] (\a,2.5)-- ++(8,0);
    \draw[dashed,line width=0.1mm] (\a,4.5)-- ++(8,0);
    \draw[dashed,line width=0.1mm] (\a,-1.5)-- ++(8,0);
    \draw[dashed,line width=0.1mm] (\a,-3.5)-- ++(8,0);
    \draw (\a,-3.5) -- (\a,4.5);
    \draw[dashed,line width=0.1mm] (2.5,-3.5) --(2.5,4.5);
    \draw[dashed,line width=0.1mm] (4.5,-3.5) --(4.5,4.5);
    \draw[dashed,line width=0.1mm] (6.5,-3.5) --(6.5,4.5);
    \draw[dashed,line width=0.1mm] (8.5,-3.5) --(8.5,4.5);
    \draw[solid,line width=0.4mm] (\a,\a)-- ++\Ed -- ++\Eu;
    \filldraw (\a,\a) circle (8pt);
    \filldraw (\a,\a)++\Ed circle (8pt);
    \filldraw (\a,\a)++\Ed++\Eu circle (8pt);
\end{tikzpicture}
\quad
\begin{tikzpicture}[scale=0.15]
    \draw[line width=0.1mm] (\a,\a)-- ++(8,0);
    \draw[dashed,line width=0.1mm] (\a,2.5)-- ++(8,0);
    \draw[dashed,line width=0.1mm] (\a,4.5)-- ++(8,0);
    \draw[dashed,line width=0.1mm] (\a,-1.5)-- ++(8,0);
    \draw[dashed,line width=0.1mm] (\a,-3.5)-- ++(8,0);
    \draw (\a,-3.5) -- (\a,4.5);
    \draw[dashed,line width=0.1mm] (2.5,-3.5) --(2.5,4.5);
    \draw[dashed,line width=0.1mm] (4.5,-3.5) --(4.5,4.5);
    \draw[dashed,line width=0.1mm] (6.5,-3.5) --(6.5,4.5);
    \draw[dashed,line width=0.1mm] (8.5,-3.5) --(8.5,4.5);
    \draw[solid,line width=0.4mm] (\a,\a)-- ++\Nd -- ++\Nu -- ++ \Nu -- ++ \Nd;
     \filldraw (\a,\a) circle (8pt);
    \filldraw (\a,\a)++\Nd circle (8pt);
    \filldraw (\a,\a)++\Nd++\Nu circle (8pt);
    \filldraw (\a,\a)++\Nd++\Nu++\Nu circle (8pt);
    \filldraw (\a,\a)++\Nd++\Nu++\Nu++\Nd circle (8pt);
\end{tikzpicture}\\\vskip0.4cm
\begin{tikzpicture}[scale=0.15]
    \draw[line width=0.1mm] (\a,\a)-- ++(8,0);
    \draw[dashed,line width=0.1mm] (\a,2.5)-- ++(8,0);
    \draw[dashed,line width=0.1mm] (\a,4.5)-- ++(8,0);
    \draw[dashed,line width=0.1mm] (\a,6.5)-- ++(8,0);
    \draw[dashed,line width=0.1mm] (\a,-1.5)-- ++(8,0);
    \draw (\a,-1.5) -- (\a,6.5);
    \draw[dashed,line width=0.1mm] (2.5,-1.5) --(2.5,6.5);
    \draw[dashed,line width=0.1mm] (4.5,-1.5) --(4.5,6.5);
    \draw[dashed,line width=0.1mm] (6.5,-1.5) --(6.5,6.5);
    \draw[dashed,line width=0.1mm] (8.5,-1.5) --(8.5,6.5);
    \draw[solid,line width=0.4mm] (\a,\a)-- ++\Nu -- ++\Nd -- ++\Nu -- ++\Nd;
    \filldraw (\a,\a) circle (8pt);
    \filldraw (\a,\a)++\Nu circle (8pt);
    \filldraw (\a,\a)++\Nu++\Nd circle (8pt);
    \filldraw (\a,\a)++\Nu++\Nd++\Nu circle (8pt);
    \filldraw (\a,\a)++\Nu++\Nd++\Nu++\Nd circle (8pt);
\end{tikzpicture}
\quad
\begin{tikzpicture}[scale=0.15]
    \draw[line width=0.1mm] (\a,\a)-- ++(8,0);
    \draw[dashed,line width=0.1mm] (\a,2.5)-- ++(8,0);
    \draw[dashed,line width=0.1mm] (\a,4.5)-- ++(8,0);
    \draw[dashed,line width=0.1mm] (\a,6.5)-- ++(8,0);
    \draw[dashed,line width=0.1mm] (\a,8.5)-- ++(8,0);
    \draw (\a,\a) -- (\a,8.5);
    \draw[dashed,line width=0.1mm] (2.5,\a) --(2.5,8.5);
    \draw[dashed,line width=0.1mm] (4.5,\a) --(4.5,8.5);
    \draw[dashed,line width=0.1mm] (6.5,\a) --(6.5,8.5);
    \draw[dashed,line width=0.1mm] (8.5,\a) --(8.5,8.5);
    \draw[solid,line width=0.4mm] (\a,\a)-- ++\Nu -- ++\Nu -- ++ \Nd -- ++ \Nd;
     \filldraw (\a,\a) circle (8pt);
    \filldraw (\a,\a)++\Nu circle (8pt);
    \filldraw (\a,\a)++\Nu++\Nu circle (8pt);
    \filldraw (\a,\a)++\Nu++\Nu++\Nd circle (8pt);
    \filldraw (\a,\a)++\Nu++\Nu++\Nd++\Nd circle (8pt);
\end{tikzpicture}
\quad
\begin{tikzpicture}[scale=0.15]
    \draw[line width=0.1mm] (\a,\a)-- ++(8,0);
    \draw[dashed,line width=0.1mm] (\a,2.5)-- ++(8,0);
    \draw[dashed,line width=0.1mm] (\a,4.5)-- ++(8,0);
    \draw[dashed,line width=0.1mm] (\a,-1.5)-- ++(8,0);
    \draw[dashed,line width=0.1mm] (\a,-3.5)-- ++(8,0);
    \draw (\a,-3.5) -- (\a,4.5);
    \draw[dashed,line width=0.1mm] (2.5,-3.5) --(2.5,4.5);
    \draw[dashed,line width=0.1mm] (4.5,-3.5) --(4.5,4.5);
    \draw[dashed,line width=0.1mm] (6.5,-3.5) --(6.5,4.5);
    \draw[dashed,line width=0.1mm] (8.5,-3.5) --(8.5,4.5);
    \draw[solid,line width=0.4mm] (\a,\a)-- ++\Eu -- ++\Ed;
     \filldraw (\a,\a) circle (8pt);
    \filldraw (\a,\a)++\Eu circle (8pt);
    \filldraw (\a,\a)++\Eu++\Ed circle (8pt);
\end{tikzpicture}
\quad
\begin{tikzpicture}[scale=0.15]
    \draw[line width=0.1mm] (\a,\a)-- ++(8,0);
    \draw[dashed,line width=0.1mm] (\a,2.5)-- ++(8,0);
    \draw[dashed,line width=0.1mm] (\a,4.5)-- ++(8,0);
    \draw[dashed,line width=0.1mm] (\a,-1.5)-- ++(8,0);
    \draw[dashed,line width=0.1mm] (\a,-3.5)-- ++(8,0);
    \draw (\a,-3.5) -- (\a,4.5);
    \draw[dashed,line width=0.1mm] (2.5,-3.5) --(2.5,4.5);
    \draw[dashed,line width=0.1mm] (4.5,-3.5) --(4.5,4.5);
    \draw[dashed,line width=0.1mm] (6.5,-3.5) --(6.5,4.5);
    \draw[dashed,line width=0.1mm] (8.5,-3.5) --(8.5,4.5);
    \draw[solid,line width=0.4mm] (\a,\a)-- ++\Nu -- ++\Nd -- ++ \Nd -- ++ \Nu;
       \filldraw (\a,\a) circle (8pt);
    \filldraw (\a,\a)++\Nu circle (8pt);
    \filldraw (\a,\a)++\Nu++\Nd circle (8pt);
    \filldraw (\a,\a)++\Nu++\Nd++\Nd circle (8pt);
    \filldraw (\a,\a)++\Nu++\Nd++\Nd++\Nu circle (8pt);
\end{tikzpicture}
\end{center}
\caption{The eight grand knight's paths of size 4 ending on the $x$-axis.}
    \label{grand knight's paths of size 4}
\end{figure}

 We refer to Figure~\ref{grand knight's paths of size 4} for the illustration of the eight grand knight's paths of size 4 ending on the $x$-axis (see the entry indicated by the boldface in the table).

\begin{corollary} The generating function for the total number $t_n$ of grand knight's paths ending at a non-negative altitude with respect to the size is 
$$H(1)=\frac{h_0u_1u_2-h_1}{(1-u_1)(1-u_2)}.$$
An asymptotic approximation of $t_n$ is 
$$\frac{\sqrt{3}}{12}\left(1+\sqrt{3}\right)^{n+1}.$$
\end{corollary} 
The first terms of $t_n$, $0\leq n\leq 12$, are 
$$1,\,1, \,4, \,8, \, 26, \,63, \,186, \,478, \,1352, \,3574, \, 9927, \,26640, \, 73354, 
$$ and this sequence does not appear in the OEIS.

\begin{corollary} The generating function for the total sum $s_n$ of the altitudes in all grand knight's paths of size $n$ ending at a non-negative altitude is  
$$\frac{\partial}{\partial u}H(u)\mid_{u=1}=\frac{h_0u_1^2u_2 + u_1u_2(h_0u_2 - 2h_0 - h_1) + h_1}{(1-u_1)^2(1-u_2)^2}.$$
An asymptotic approximation of $s_n$ is
$$\frac{(4\sqrt{3}+7)\sqrt{2}\sqrt{137\sqrt{3}-237}}{6}\sqrt{\frac{n}{\pi}}\left(1+\sqrt{3}\right)^n.$$
\end{corollary}
The first terms of $s_n$, $0\leq n\leq 12$, are 
$$0,\, 2,\,5, \,20, \,56,\, 180, \,516, \,1552, \,4452, \,13000, \,37120, \,106684, \,303090, $$ and this sequence does not appear in OEIS.

\begin{corollary}
    The expected altitude of a grand knight's path of size $n$ ending at positive altitude is asymptotically
    $$(5+3\sqrt{3})\sqrt{137\sqrt{3}-237}\sqrt{\frac{2n}{3\pi}}.$$
\end{corollary}

\section{Unrestricted grand zigzag knight's paths}\label{unrestricted zigzag}

Recall from Definition~\ref{d12} that a grand zigzag knight's path is a grand knight's path avoiding the consecutive patterns $NN$, $NE$, $\bar{N}\bar{N}$, $\bar{N}\bar{E}$, $EE$, $EN$, $\bar{E}\bar{E}$, $\bar{E}\bar{N}$. 

Let $\mathcal{Z}_{n,k}$ be the set of grand zigzag knight's paths of size $n$ and altitude $k$, and $\mathcal{Z}_{n,k}^+$ (resp.~$\mathcal{Z}_{n,k}^-)$ be the subset of $\mathcal{Z}_{n,k}$ of paths starting with $E$ or $N$ (resp.~$\bar{E}$ or $\bar{N}$). Note that the symmetry with respect to the $x$-axis provides a direct bijection between $\mathcal{Z}_{n,k}$ and $\mathcal{Z}_{n,-k}$, and  also between $\mathcal{Z}_{n,k}^+$ and $\mathcal{Z}_{n,-k}^-$.

\subsection{Enumeration using algebraic method}\label{algmeth}
Let $\hat{F}(x,u,z)$ be the generating
function for the number of grand zigzag knight's paths ending with $E$ or $N$ with respect to the number of steps (marked by $x$), the altitude (marked by $u$, the power of $u$ is negative for paths ending below $x$-axis) and the size (final $x$-coordinate, marked by $z$). Denote by $\hat{G}(x,u,z)$ the generating function for grand zigzag knight's paths ending with $\bar{E}$ or $\bar{N}$ and by $\hat{Z}(x,u,z)$ the generating function for all grand zigzag knight paths.

Then we have readily the following system of equations:
\begin{equation}
    \begin{cases}
    \hat{F}(x,u,z) & = (1 + \hat{G}(x,u,z)) 
    (x z u^2 + x z^2 u ),\\
    \hat{G}(x,u,z) & = (1 + \hat{F}(x,u,z)) 
    (x z u^{-2} + x z^2 u^{-1}),\\
    \hat{Z}(x,u,z) & = 1 + \hat{F}(x,u,z) + \hat{G}(x,u,z).
    \end{cases}\label{Sys}
\end{equation}
The solutions of (\ref{Sys}) are rational:
\begin{equation*}
\begin{aligned}
\hat{F}(x,u,z) & = - \frac{x z \left(u + z\right) \left(u^{2} + u x z^{2} + x z\right)}{u^{2} x^{2} z^{3} + u x^{2} z^{4} + u x^{2} z^{2} - u + x^{2} z^{3}},\\
\hat{G}(x,u,z) & = - \frac{x z \left(u z + 1\right) \left(u^{2} x z + u x z^{2} + 1\right)}{u \left(u^{2} x^{2} z^{3} + u x^{2} z^{4} + u x^{2} z^{2} - u + x^{2} z^{3}\right)},\\
\hat{Z}(x,u,z) & = - \frac{\left(u^{2} + u x z^{2} + x z\right) \left(u^{2} x z + u x z^{2} + 1\right)}{u \left(u^{2} x^{2} z^{3} + u x^{2} z^{4} + u x^{2} z^{2} - u + x^{2} z^{3}\right
)}.
\end{aligned}
\end{equation*}

Fixing $u=1$ and $x=1$ in the third equation, we obtain the generating function for all grand zigzag knight's paths of a given size:
\begin{equation}\label{zsec3}\hat{Z}(1,1,z)  = \frac{1 +z+ z^{2} }{1-z-z^2}.
\end{equation}
 Here are the first coefficients of $z^n$ for $0\leq n\leq 16$:
 $$1,\, 2, \, 4, \, 6, \, 10, \, 16, \, 26, \, 42, \, 68, \, 110, \, 178, \, 288, \, 466, \, 754, \, 1220, \, 1974, \, 3194.$$
This corresponds to the sequence \href{https://oeis.org/A128588}{A128588} in \cite{OEIS}, which is a Fibonacci-like sequence with $a_0=1$, $a_1=2$ and $a_2=4$ (the $n$-th term, except the first, is twice the Fibonacci number $F_n$ defined by $F_n=F_{n-1}+F_{n-2}$ for $n\geq 2$, with $F_0=0$ and $F_1=1$).

By decomposing the rational fraction $\hat{Z}(x,u,z)$
into partial fractions, and after isolating the terms which are analytic 
 at $u=0$, we can obtain an expression for the 
generating function for paths of positive height. In this paper, we use a different technique inspired by Prodinger~\cite{Pro} and based on the kernel method.
 We describe it in full details
in order to introduce the routine that we will use in Section~4 for bounded paths.

For $k\in\Z$, let $f_k$ (resp.~$g_k$) be the generating function of the number of grand zigzag knight's paths of altitude $k$ with $E$ or $N$ (resp.~$\bar{E}$ or $\bar{N}$), with respect to the size. Due to the symmetry with respect to the $x$-axis, we have $f_{-k}=g_k$ for $k\ne 0$, and $f_0=1+g_0$. So, we will focus on positive altitudes. Let us consider $F(u,z)=\sum_{k\geq 0}f_k(z)u^k$ and $G(u,z)=\sum_{k\geq 0}g_k(z)u^k$ (for short, we will write $F(u)$ and $G(u)$). We easily obtain the following equations, with the convention that the empty path is counted in $f_0$:
\begin{align}
f_1&=z^2+zg_{-1}+z^2g_0=z^2+zf_1+z^2g_0=zf_1+z^2f_0,\label{f_1}\\ 
f_2&=z+zg_0+z^2g_1, \label{f_2}\\
f_k&=zg_{k-2}+z^2g_{k-1} \  \text{ for } k\geq3, \label{f_k}\\
g_k&=zf_{k+2}+z^2f_{k+1} \ \text{ for } k\geq0. \label{g_k}
\end{align}
From equations~(\ref{f_1}-\ref{g_k}) we deduce the following functional equations:
\begin{align*}
    F(u)&=\left(1+\frac{z^3u}{1-z}\right)f_0+zu(z+u)(G(u)+1),\\
    G(u)&=-\left(\frac{z}{u^2}+\frac{z^2}{u}+\frac{z^3}{u(1-z)}\right)f_0+\left(\frac{z}{u^2}+\frac{z^2}{u}\right)F(u),
\end{align*}
which yields 
\begin{align*}
F(u)&= \scalemath{0.85}{{\frac {f_0(-\,u{z}^{4}+\,u{z}^{3}+\,{z}^{4} -\,u{z}^{2} -\,{z}^{3} -\,uz +\,u)   -{u}^{3}{z}^{2
}-{u}^{2}{z}^{3}+{u}^{3}z+{u}^{2}
{z}^{2}
}{
 \left( -1+z \right)  \left( {u}^{2}{z}^{3}+u{z}^{4}+{z}^{2}u+{z}^{3}-
u \right) }}},\\
G(u)&= \scalemath{0.85}{{\frac { \left( f_0(\,u{z}^{3}+\,{
z}^{2}  -z       )-{u}^{2}{z}^{2}-{z}^{3}u+{u}^{2}z+{z}^{2}u-zu-{z}^{2}+u+z \right) {z}^{2}}{
 \left( -1+z \right)  \left( {u}^{2}{z}^{3}+u{z}^{4}+{z}^{2}u+{z}^{3}-
u \right) }}}.
\end{align*}
Let
\begin{align*}
&r:={\frac {1-{z}^{4}-{z}^{2}-\sqrt {{z}^{8}-2\,{z}^{6}-{z}^{4}-2\,{
z}^{2}+1}}{2{z}^{3}}}, \\
&s:={\frac {1-{z}^{4}-{z}^{2}+\sqrt {{z}^{8}-2\,{z}^{6}-{z}^{4}-2\,{z}
^{2}+1}}{2{z}^{3}}}
\end{align*}
be the two roots in $u$ of the kernel $ {u}^{2}{z}^{3}+u{z}^{4}+{z}^{2}u+{z}^{3}-
u$.
Applying the kernel method, we obtain
\begin{align*}f_0&=\frac{r(z-1)}{z^3(rz^2+z-1)}
=\frac{(1 - z)(-1 + z^2 + z^4 +\sqrt{1 - 2 z^2 - z^4 - 2 z^6 + z^8})}{z^5(1 - 2 z + z^2 - z^4 -\sqrt{1 - 2 z^2 - z^4 - 2 z^6 + z^8})},
\end{align*}
and we can state the following cancelling the factor $(u-r)$.

\begin{theorem}
The generating functions $F(u)$, $G(u)$ are given by:
\begin{align*}F(u)&=\frac{-u^2-u(r+z)-z^2f_0s}{z^2(u-s)} \quad \text{and} \quad
G(u)=\frac{-z^2u-z^2f_0s+z^2s}{z^2(u-s)}.
\end{align*}
The bivariate generating function $H(u):=F(u)+G(u)$ for the number of grand zigzag knight's paths (ending at a non-negative $y$-coordinate) with respect to the size and the altitude is
$$H(u)=-\frac{u^2+u(r+z+z^2)+z^2s(2f_0-1)}{z^2(u-s)}.$$
The generating function for the total number of grand zigzag knight's paths of size $n$ ending at a non negative $y$-coordinate is $$H(1)=-\frac{1+r+z+z^2+z^2s(2f_0-1)}{z^2(1-s)}.$$ 
\end{theorem}
Here are the first terms of $H(1)$ for $0\leq n\leq 16$: 
    $$1,\, 1, \, 3, \, 3, \, 7, \, 9, \, 18, \, 24, \, 45, \, 63, \, 115, \, 166, \, 296, \, 435, \, 763, \, 1138, \, 1973.$$
    Note that the generating function for the total number of grand zigzag knight's paths is thus $$2H(1)-(2f_0-1)=\frac{1 +z+ z^{2} }{1-z-z^2}.$$ Consistently with the beginning of this section, the coefficient of $z^n$ in $2H(1)-(2f_0-1)$ is twice the $n$-th Fibonacci number (see \href{https://oeis.org/A128588}{A128588} in \cite{OEIS}). This can be also seen bijectively as follows. A grand zigzag knight's paths of size $n$ starting with an up-step (so half of the total number of grand zigzag knight's paths of size $n$) can start either with $E$, and then followed by a grand zigzag knight's path of size $n-2$ starting with a down-step, or it can start with $N$, and then followed by a grand zigzag knight's path of size $n-1$ starting with a down-step. See Table~\ref{values Znk} for the first  values of $|\mathcal{Z}_{n,k}|$. Note that the sequence $(|\mathcal{Z}_{2n,0}^+|)_{n\geq0}$ corresponds to \href{https://oeis.org/A051286}{A051286} in~\cite{OEIS}.

\begin{corollary}\label{zigzag ending at height k}
The generating function for the number of grand zigzag knight's paths ending at $y$-coordinate $k$ with respect to the size is given by:
\begin{align*}[u^0]H(u)&=2f_0-1,\quad [u^1]H(u)=\frac{r+z+2z^2f_0}{z^2s},\\
[u^k]H(u)&=\frac{r^{k-1}}{z^2}(1+r(r+z)+2rz^2f_0)\quad\mbox{for } k\geq 2.
\end{align*}
\end{corollary}

\begin{corollary} An asymptotic approximation for the expected altitude of a grand zigzag knight's path of size $n$ ending on or above the $x$-axis is 
$$\frac{2(\sqrt{5}-2)}{\sqrt{7\sqrt{5}-15}}\sqrt{\frac{n}{\pi}}.$$
 This is also the variance (divided by 2) of the altitude  of a grand zigzag knight's path (not necessarily ending on or above the $x$-axis) of size $n$ (the expected value of the altitude  being simply $0$ for those paths).
\end{corollary}
\begin{proof}
    The generating function for the total number of grand zigzag knight's paths of size $n$ is $H(1)$, and the one for the total sum of altitudes over all grand zigzag knight's paths of size $n$ is $\frac{\partial}{\partial u}H(u)\mid_{u=1}$. We compute asymptotics of the coefficients of both those generating functions using singularity analysis (see \cite{analytic combinatorics}), and then we take the ratio to obtain the stated result.
\end{proof}

\begin{rmk}
    From \cite[Theorem 1]{BarRa}, we deduce that the expected altitude  of a grand zigzag knight's path of size $n$ and staying above the $x$-axis is equivalent as $n\to\infty$ to
    $$\frac{(5+\sqrt{5})\sqrt{7\sqrt{5}-15}}{20}\sqrt{\pi n}.$$
    Thus, the altitude of a zigzag knight's path (i.e. a grand zigzag knight's path in $\mathbb{N}^2$), is, in average, $\approx 1.57079633$ times the altitude of a grand zigzag knight's path in $\mathbb{Z}^2$ ending at non-negative altitude.
\end{rmk}

\begin{table}[ht]
\centering
\rowcolors{2}{}{lightgray}
\begin{tabular}{|c|cccccccccccccccc|}
    \hline
    \diagbox[]{$k$}{$n$} & 0 & 1 & 2 & 3 & 4 & 5 & 6 & 7 & 8 & 9 & 10 & 11 & 12 & 13 & 14 & 15\\
    \hline
    0 & 1 & 0 & 2 & 0 & 4 & 2 & 10 & $\mathbf{6}$ & 22 & 16 & 52 & 44 & 126 & 116 & 306 & 302\\
    1 & 0 & 0 & 1 & 2 & 2 & 4 & 4 & 10 & 11 & 26 & 28 & 64 & 71 & 160 & 183 & 402\\
    2 & 0 & 1 & 0 & 1 & 0 & 3 & 2 & 7 & 6 & 16 & 18 & 40 & 52 & 100 & 142 & 252\\
    3 & 0 & 0 & 0 & 0 & 1 & 0 & 2 & 0 & 6 & 2 & 16 & 8 & 41 & 28 & 107 & 90\\
    4 & 0 & 0 & 0 & 0 & 0 & 0 & 0 & 1 & 0 & 3 & 0 & 10 & 2 & 30 & 10 & 85\\
    \hline
\end{tabular}
\caption{The number of grand zigzag knight's paths from $(0,0)$ to $(n,k)$ for $(n,k)\in[0,15]\times[0,4]$.}\label{values Znk}
\end{table}

\begin{figure}[ht]
\centering
     \begin{tikzpicture}[scale=0.15]
    \draw (\a,\a)-- ++(14,0);
    \draw[dashed,line width=0.1mm] (\a,2.5)-- ++(14,0);
    \draw[dashed,line width=0.1mm] (\a,4.5)-- ++(14,0);
    \draw[dashed,line width=0.1mm] (\a,-1.5)-- ++(14,0);
    \draw[dashed,line width=0.1mm] (\a,-3.5)-- ++(14,0);
    \draw (\a,-3.5) -- (\a,4.5);
    \draw[dashed,line width=0.1mm] (2.5,-3.5) --(2.5,4.5);
    \draw[dashed,line width=0.1mm] (4.5,-3.5) --(4.5,4.5);
    \draw[dashed,line width=0.1mm] (6.5,-3.5) --(6.5,4.5);
    \draw[dashed,line width=0.1mm] (8.5,-3.5) --(8.5,4.5);
    \draw[dashed,line width=0.1mm] (10.5,-3.5) --(10.5,4.5);
    \draw[dashed,line width=0.1mm] (12.5,-3.5) --(12.5,4.5);
    \draw[dashed,line width=0.1mm] (14.5,-3.5) --(14.5,4.5);
    \draw[solid,line width=0.4mm] (\a,\a)-- ++\Eu--++\Nd-- ++\Nu-- ++\Nd-- ++\Eu;
     \filldraw (\a,\a) circle (8pt);
    \filldraw (\a,\a)++\Eu circle (8pt);
    \filldraw (\a,\a)++\Eu++\Nd circle (8pt);
    \filldraw (\a,\a)++\Eu++\Nd++\Nu circle (8pt);
    \filldraw (\a,\a)++\Eu++\Nd++\Nu++\Nd circle (8pt);
    \filldraw (\a,\a)++\Eu++\Nd++\Nu++\Nd++\Eu circle (8pt);
\end{tikzpicture}
\quad
\begin{tikzpicture}[scale=0.15]
    \draw (\a,\a)-- ++(14,0);
    \draw[dashed,line width=0.1mm] (\a,2.5)-- ++(14,0);
    \draw[dashed,line width=0.1mm] (\a,4.5)-- ++(14,0);
    \draw[dashed,line width=0.1mm] (\a,-1.5)-- ++(14,0);
    \draw[dashed,line width=0.1mm] (\a,-3.5)-- ++(14,0);
    \draw (\a,-3.5) -- (\a,4.5);
    \draw[dashed,line width=0.1mm] (2.5,-3.5) --(2.5,4.5);
    \draw[dashed,line width=0.1mm] (4.5,-3.5) --(4.5,4.5);
    \draw[dashed,line width=0.1mm] (6.5,-3.5) --(6.5,4.5);
    \draw[dashed,line width=0.1mm] (8.5,-3.5) --(8.5,4.5);
    \draw[dashed,line width=0.1mm] (10.5,-3.5) --(10.5,4.5);
    \draw[dashed,line width=0.1mm] (12.5,-3.5) --(12.5,4.5);
    \draw[dashed,line width=0.1mm] (14.5,-3.5) --(14.5,4.5);
    \draw[solid,line width=0.4mm] (\a,\a)-- ++\Eu-- ++\Nd-- ++\Eu -- ++\Nd--++\Nu;
    \filldraw (\a,\a) circle (8pt);
    \filldraw (\a,\a)++\Eu circle (8pt);
    \filldraw (\a,\a)++\Eu++\Nd circle (8pt);
    \filldraw (\a,\a)++\Eu++\Nd++\Eu circle (8pt);
    \filldraw (\a,\a)++\Eu++\Nd++\Eu++\Nd circle (8pt);
    \filldraw (\a,\a)++\Eu++\Nd++\Eu++\Nd++\Nu circle (8pt);
\end{tikzpicture}
\quad
\begin{tikzpicture}[scale=0.15]
    \draw (\a,\a)-- ++(14,0);
    \draw[dashed,line width=0.1mm] (\a,2.5)-- ++(14,0);
    \draw[dashed,line width=0.1mm] (\a,4.5)-- ++(14,0);
    \draw[dashed,line width=0.1mm] (\a,-1.5)-- ++(14,0);
    \draw[dashed,line width=0.1mm] (\a,-3.5)-- ++(14,0);
    \draw (\a,-3.5) -- (\a,4.5);
    \draw[dashed,line width=0.1mm] (2.5,-3.5) --(2.5,4.5);
    \draw[dashed,line width=0.1mm] (4.5,-3.5) --(4.5,4.5);
    \draw[dashed,line width=0.1mm] (6.5,-3.5) --(6.5,4.5);
    \draw[dashed,line width=0.1mm] (8.5,-3.5) --(8.5,4.5);
    \draw[dashed,line width=0.1mm] (10.5,-3.5) --(10.5,4.5);
    \draw[dashed,line width=0.1mm] (12.5,-3.5) --(12.5,4.5);
    \draw[dashed,line width=0.1mm] (14.5,-3.5) --(14.5,4.5);
    \draw[solid,line width=0.4mm] (\a,\a)-- ++\Nu--++\Nd-- ++\Eu-- ++\Nd-- ++\Eu;
    \filldraw (\a,\a) circle (8pt);
    \filldraw (\a,\a)++\Nu circle (8pt);
    \filldraw (\a,\a)++\Nu++\Nd circle (8pt);
    \filldraw (\a,\a)++\Nu++\Nd++\Eu circle (8pt);
    \filldraw (\a,\a)++\Nu++\Nd++\Eu++\Nd circle (8pt);
    \filldraw (\a,\a)++\Nu++\Nd++\Eu++\Nd++\Eu circle (8pt);
\end{tikzpicture}\\\vskip0.4cm
\begin{tikzpicture}[scale=0.15]
    \draw (\a,\a)-- ++(14,0);
    \draw[dashed,line width=0.1mm] (\a,2.5)-- ++(14,0);
    \draw[dashed,line width=0.1mm] (\a,4.5)-- ++(14,0);
    \draw[dashed,line width=0.1mm] (\a,-1.5)-- ++(14,0);
    \draw[dashed,line width=0.1mm] (\a,-3.5)-- ++(14,0);
    \draw (\a,-3.5) -- (\a,4.5);
    \draw[dashed,line width=0.1mm] (2.5,-3.5) --(2.5,4.5);
    \draw[dashed,line width=0.1mm] (4.5,-3.5) --(4.5,4.5);
    \draw[dashed,line width=0.1mm] (6.5,-3.5) --(6.5,4.5);
    \draw[dashed,line width=0.1mm] (8.5,-3.5) --(8.5,4.5);
    \draw[dashed,line width=0.1mm] (10.5,-3.5) --(10.5,4.5);
    \draw[dashed,line width=0.1mm] (12.5,-3.5) --(12.5,4.5);
    \draw[dashed,line width=0.1mm] (14.5,-3.5) --(14.5,4.5);
    \draw[solid,line width=0.4mm] (\a,\a)-- ++\Ed--++\Nu-- ++\Nd-- ++\Nu-- ++\Ed;
     \filldraw (\a,\a) circle (8pt);
    \filldraw (\a,\a)++\Ed circle (8pt);
    \filldraw (\a,\a)++\Ed++\Nu circle (8pt);
    \filldraw (\a,\a)++\Ed++\Nu++\Nd circle (8pt);
    \filldraw (\a,\a)++\Ed++\Nu++\Nd++\Nu circle (8pt);
     \filldraw (\a,\a)++\Ed++\Nu++\Nd++\Nu++\Ed circle (8pt);
\end{tikzpicture}
\quad
\begin{tikzpicture}[scale=0.15]
    \draw (\a,\a)-- ++(14,0);
    \draw[dashed,line width=0.1mm] (\a,2.5)-- ++(14,0);
    \draw[dashed,line width=0.1mm] (\a,4.5)-- ++(14,0);
    \draw[dashed,line width=0.1mm] (\a,-1.5)-- ++(14,0);
    \draw[dashed,line width=0.1mm] (\a,-3.5)-- ++(14,0);
    \draw (\a,-3.5) -- (\a,4.5);
    \draw[dashed,line width=0.1mm] (2.5,-3.5) --(2.5,4.5);
    \draw[dashed,line width=0.1mm] (4.5,-3.5) --(4.5,4.5);
    \draw[dashed,line width=0.1mm] (6.5,-3.5) --(6.5,4.5);
    \draw[dashed,line width=0.1mm] (8.5,-3.5) --(8.5,4.5);
    \draw[dashed,line width=0.1mm] (10.5,-3.5) --(10.5,4.5);
    \draw[dashed,line width=0.1mm] (12.5,-3.5) --(12.5,4.5);
    \draw[dashed,line width=0.1mm] (14.5,-3.5) --(14.5,4.5);
    \draw[solid,line width=0.4mm] (\a,\a)-- ++\Ed-- ++\Nu-- ++\Ed -- ++\Nu--++\Nd;
    \filldraw (\a,\a) circle (8pt);
    \filldraw (\a,\a)++\Ed circle (8pt);
    \filldraw (\a,\a)++\Ed++\Nu circle (8pt);
    \filldraw (\a,\a)++\Ed++\Nu++\Ed circle (8pt);
    \filldraw (\a,\a)++\Ed++\Nu++\Ed++\Nu circle (8pt);
    \filldraw (\a,\a)++\Ed++\Nu++\Ed++\Nu++\Nd circle (8pt);
\end{tikzpicture}
\quad
\begin{tikzpicture}[scale=0.15]
    \draw (\a,\a)-- ++(14,0);
    \draw[dashed,line width=0.1mm] (\a,2.5)-- ++(14,0);
    \draw[dashed,line width=0.1mm] (\a,4.5)-- ++(14,0);
    \draw[dashed,line width=0.1mm] (\a,-1.5)-- ++(14,0);
    \draw[dashed,line width=0.1mm] (\a,-3.5)-- ++(14,0);
    \draw (\a,-3.5) -- (\a,4.5);
    \draw[dashed,line width=0.1mm] (2.5,-3.5) --(2.5,4.5);
    \draw[dashed,line width=0.1mm] (4.5,-3.5) --(4.5,4.5);
    \draw[dashed,line width=0.1mm] (6.5,-3.5) --(6.5,4.5);
    \draw[dashed,line width=0.1mm] (8.5,-3.5) --(8.5,4.5);
    \draw[dashed,line width=0.1mm] (10.5,-3.5) --(10.5,4.5);
    \draw[dashed,line width=0.1mm] (12.5,-3.5) --(12.5,4.5);
    \draw[dashed,line width=0.1mm] (14.5,-3.5) --(14.5,4.5);
    \draw[solid,line width=0.4mm] (\a,\a)-- ++\Nd--++\Nu-- ++\Ed-- ++\Nu-- ++\Ed;
    \filldraw (\a,\a) circle (8pt);
    \filldraw (\a,\a)++\Nd circle (8pt);
    \filldraw (\a,\a)++\Nd++\Nu circle (8pt);
    \filldraw (\a,\a)++\Nd++\Nu++\Ed circle (8pt);
    \filldraw (\a,\a)++\Nd++\Nu++\Ed++\Nu circle (8pt);
    \filldraw (\a,\a)++\Nd++\Nu++\Ed++\Nu++\Ed circle (8pt);
\end{tikzpicture}
    \caption{The six grand zigzag knight's paths of size 7 ending on the $x$-axis.}
    \label{zigzag of size 7 and final height 0}
\end{figure}
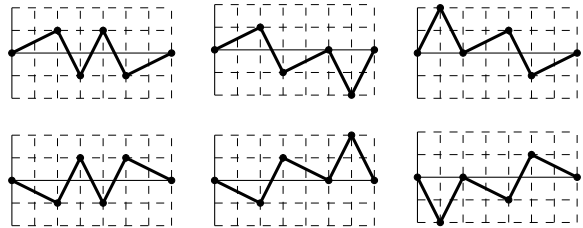

We refer to Figure~\ref{zigzag of size 7 and final height 0} for the illustration of the six grand zigzag knight's paths of size $7$ ending on the $x$-axis (see the entry indicated by the boldface in the table).

Now, we focus on grand zigzag knight's paths that never touch the $x$-axis except for the initial and the final point. Let $\mathcal{B}$ be the set of such paths.

\begin{corollary}
    The generating function $B(z)$ for the number $b_n$ of grand zigzag knight's paths in $\mathcal{B}$ with respect to the size is 
    $$B(z)=\frac{2z^5r+2z^4-2z^3-3rz+3r}{r(1-z)}.$$
\end{corollary}
\begin{proof}
    Let $\Gamma=2f_0-1$ be the generating function for the grand zigzag knight's paths ending on the $x$-axis. Each element counted by $\Gamma$ is either empty or can be uniquely decomposed into a juxtaposition $B_1\cdots B_k$ of non-empty paths $B_i\in\mathcal{B}$ so that $B_iB_{i+1}$ forms a zigzag, which forces $B_{i+1}$ to start with an up step (resp.~down step) if $B_i$ ends with a down step (resp.~up step). Therefore, the generating function for $B_1$ is $B(z)-1$, and the generating function for $B_i$, $i\geq 2$, is $(B(z)-1)/2$. Thus, we have 
   $$\Gamma=1+(B(z)-1)\sum_{k\geq 0}\left(\frac{B(z)-1}{2}\right)^k=\frac{1+B(z)}{3-B(z)}.$$
   We deduce $B(z)=(3\Gamma-1)/(\Gamma+1)$ and the expression of $B(z)$ follows from Corollary~\ref{zigzag ending at height k}.
\end{proof}

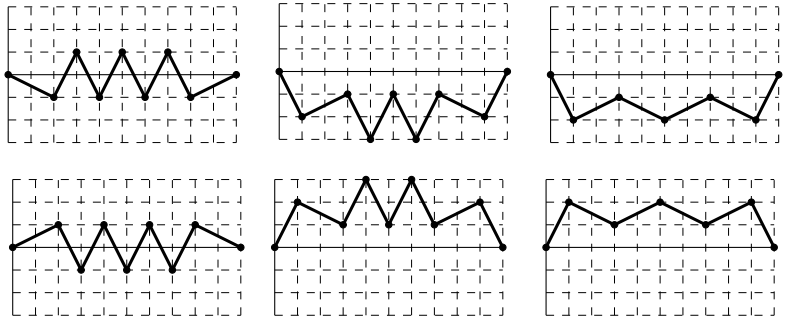
\begin{figure}[ht]
\centering
    \begin{center}
        \begin{tikzpicture}[scale=0.15]
    \draw (\a,\a)-- ++(20,0);
    \draw[dashed,line width=0.1mm] (\a,2.5)-- ++(20,0);
    \draw[dashed,line width=0.1mm] (\a,4.5)-- ++(20,0);
    \draw[dashed,line width=0.1mm] (\a,6.5)-- ++(20,0);
    \draw[dashed,line width=0.1mm] (\a,-1.5)-- ++(20,0);
    \draw[dashed,line width=0.1mm] (\a,-3.5)-- ++(20,0);
    \draw[dashed,line width=0.1mm] (\a,-5.5)-- ++(20,0);
    \draw (\a,-5.5) -- (\a,6.5);
    \draw[dashed,line width=0.1mm] (2.5,-5.5) --(2.5,6.5);
    \draw[dashed,line width=0.1mm] (4.5,-5.5) --(4.5,6.5);
    \draw[dashed,line width=0.1mm] (6.5,-5.5) --(6.5,6.5);
    \draw[dashed,line width=0.1mm] (8.5,-5.5) --(8.5,6.5);
    \draw[dashed,line width=0.1mm] (10.5,-5.5) --(10.5,6.5);
    \draw[dashed,line width=0.1mm] (12.5,-5.5) --(12.5,6.5);
    \draw[dashed,line width=0.1mm] (14.5,-5.5) --(14.5,6.5);
    \draw[dashed,line width=0.1mm] (16.5,-5.5) --(16.5,6.5);
    \draw[dashed,line width=0.1mm] (18.5,-5.5) --(18.5,6.5);
    \draw[dashed,line width=0.1mm] (20.5,-5.5) --(20.5,6.5);
    \draw[solid,line width=0.4mm] (\a,\a)-- ++\Ed--++\Nu-- ++\Nd-- ++\Nu-- ++\Nd-- ++\Nu-- ++\Nd-- ++\Eu;
     \filldraw (\a,\a) circle (8pt);
    \filldraw (\a,\a)++\Ed circle (8pt);
    \filldraw (\a,\a)++\Ed++\Nu circle (8pt);
    \filldraw (\a,\a)++\Ed++\Nu++\Nd circle (8pt);
    \filldraw (\a,\a)++\Ed++\Nu++\Nd++\Nu circle (8pt);
    \filldraw (\a,\a)++\Ed++\Nu++\Nd++\Nu++\Nd circle (8pt);
    \filldraw (\a,\a)++\Ed++\Nu++\Nd++\Nu++\Nd++\Nu circle (8pt);
    \filldraw (\a,\a)++\Ed++\Nu++\Nd++\Nu++\Nd++\Nu++\Nd circle (8pt);
     \filldraw (\a,\a)++\Ed++\Nu++\Nd++\Nu++\Nd++\Nu++\Nd++\Eu circle (8pt);
\end{tikzpicture}
\quad
\begin{tikzpicture}[scale=0.15]
    \draw (\a,\a)-- ++(20,0);
    \draw[dashed,line width=0.1mm] (\a,2.5)-- ++(20,0);
    \draw[dashed,line width=0.1mm] (\a,4.5)-- ++(20,0);
    \draw[dashed,line width=0.1mm] (\a,6.5)-- ++(20,0);
    \draw[dashed,line width=0.1mm] (\a,-1.5)-- ++(20,0);
    \draw[dashed,line width=0.1mm] (\a,-3.5)-- ++(20,0);
    \draw[dashed,line width=0.1mm] (\a,-5.5)-- ++(20,0);
    \draw (\a,-5.5) -- (\a,6.5);
    \draw[dashed,line width=0.1mm] (2.5,-5.5) --(2.5,6.5);
    \draw[dashed,line width=0.1mm] (4.5,-5.5) --(4.5,6.5);
    \draw[dashed,line width=0.1mm] (6.5,-5.5) --(6.5,6.5);
    \draw[dashed,line width=0.1mm] (8.5,-5.5) --(8.5,6.5);
    \draw[dashed,line width=0.1mm] (10.5,-5.5) --(10.5,6.5);
    \draw[dashed,line width=0.1mm] (12.5,-5.5) --(12.5,6.5);
    \draw[dashed,line width=0.1mm] (14.5,-5.5) --(14.5,6.5);
    \draw[dashed,line width=0.1mm] (16.5,-5.5) --(16.5,6.5);
    \draw[dashed,line width=0.1mm] (18.5,-5.5) --(18.5,6.5);
    \draw[dashed,line width=0.1mm] (20.5,-5.5) --(20.5,6.5);
    \draw[solid,line width=0.4mm] (\a,\a)-- ++\Nd-- ++\Eu -- ++\Nd -- ++\Nu -- ++\Nd -- ++\Nu -- ++\Ed -- ++\Nu;
      \filldraw (\a,\a) circle (8pt);
    \filldraw (\a,\a)++\Nd circle (8pt);
    \filldraw (\a,\a)++\Nd++\Eu circle (8pt);
    \filldraw (\a,\a)++\Nd++\Eu++\Nd circle (8pt);
    \filldraw (\a,\a)++\Nd++\Eu++\Nd++\Nu circle (8pt);
    \filldraw (\a,\a)++\Nd++\Eu++\Nd++\Nu++\Nd circle (8pt);
    \filldraw (\a,\a)++\Nd++\Eu++\Nd++\Nu++\Nd++\Nu circle (8pt);
    \filldraw (\a,\a)++\Nd++\Eu++\Nd++\Nu++\Nd++\Nu++\Ed circle (8pt);
     \filldraw (\a,\a)++\Nd++\Eu++\Nd++\Nu++\Nd++\Nu++\Ed++\Nu circle (8pt);
\end{tikzpicture}
\quad
\begin{tikzpicture}[scale=0.15]
    \draw (\a,\a)-- ++(20,0);
    \draw[dashed,line width=0.1mm] (\a,2.5)-- ++(20,0);
    \draw[dashed,line width=0.1mm] (\a,4.5)-- ++(20,0);
    \draw[dashed,line width=0.1mm] (\a,6.5)-- ++(20,0);
    \draw[dashed,line width=0.1mm] (\a,-1.5)-- ++(20,0);
    \draw[dashed,line width=0.1mm] (\a,-3.5)-- ++(20,0);
    \draw[dashed,line width=0.1mm] (\a,-5.5)-- ++(20,0);
    \draw (\a,-5.5) -- (\a,6.5);
    \draw[dashed,line width=0.1mm] (2.5,-5.5) --(2.5,6.5);
    \draw[dashed,line width=0.1mm] (4.5,-5.5) --(4.5,6.5);
    \draw[dashed,line width=0.1mm] (6.5,-5.5) --(6.5,6.5);
    \draw[dashed,line width=0.1mm] (8.5,-5.5) --(8.5,6.5);
    \draw[dashed,line width=0.1mm] (10.5,-5.5) --(10.5,6.5);
    \draw[dashed,line width=0.1mm] (12.5,-5.5) --(12.5,6.5);
    \draw[dashed,line width=0.1mm] (14.5,-5.5) --(14.5,6.5);
    \draw[dashed,line width=0.1mm] (16.5,-5.5) --(16.5,6.5);
    \draw[dashed,line width=0.1mm] (18.5,-5.5) --(18.5,6.5);
    \draw[dashed,line width=0.1mm] (20.5,-5.5) --(20.5,6.5);
    \draw[solid,line width=0.4mm] (\a,\a)-- ++\Nd-- ++\Eu -- ++\Ed-- ++\Eu -- ++\Ed -- ++\Nu;
       \filldraw (\a,\a) circle (8pt);
    \filldraw (\a,\a)++\Nd circle (8pt);
    \filldraw (\a,\a)++\Nd++\Eu circle (8pt);
    \filldraw (\a,\a)++\Nd++\Eu++\Ed circle (8pt);
    \filldraw (\a,\a)++\Nd++\Eu++\Ed++\Eu circle (8pt);
    \filldraw (\a,\a)++\Nd++\Eu++\Ed++\Eu++\Ed circle (8pt);
    \filldraw (\a,\a)++\Nd++\Eu++\Ed++\Eu++\Ed++\Nu circle (8pt);
\end{tikzpicture}\\\vskip0.4cm
 \begin{tikzpicture}[scale=0.15]
    \draw (\a,\a)-- ++(20,0);
    \draw[dashed,line width=0.1mm] (\a,2.5)-- ++(20,0);
    \draw[dashed,line width=0.1mm] (\a,4.5)-- ++(20,0);
    \draw[dashed,line width=0.1mm] (\a,6.5)-- ++(20,0);
    \draw[dashed,line width=0.1mm] (\a,-1.5)-- ++(20,0);
    \draw[dashed,line width=0.1mm] (\a,-3.5)-- ++(20,0);
    \draw[dashed,line width=0.1mm] (\a,-5.5)-- ++(20,0);
    \draw (\a,-5.5) -- (\a,6.5);
    \draw[dashed,line width=0.1mm] (2.5,-5.5) --(2.5,6.5);
    \draw[dashed,line width=0.1mm] (4.5,-5.5) --(4.5,6.5);
    \draw[dashed,line width=0.1mm] (6.5,-5.5) --(6.5,6.5);
    \draw[dashed,line width=0.1mm] (8.5,-5.5) --(8.5,6.5);
    \draw[dashed,line width=0.1mm] (10.5,-5.5) --(10.5,6.5);
    \draw[dashed,line width=0.1mm] (12.5,-5.5) --(12.5,6.5);
    \draw[dashed,line width=0.1mm] (14.5,-5.5) --(14.5,6.5);
    \draw[dashed,line width=0.1mm] (16.5,-5.5) --(16.5,6.5);
    \draw[dashed,line width=0.1mm] (18.5,-5.5) --(18.5,6.5);
    \draw[dashed,line width=0.1mm] (20.5,-5.5) --(20.5,6.5);
    \draw[solid,line width=0.4mm] (\a,\a)-- ++\Eu-- ++\Nd-- ++\Nu-- ++\Nd--++\Nu-- ++\Nd-- ++\Nu-- ++\Ed;
    \filldraw (\a,\a) circle (8pt);
    \filldraw (\a,\a)++\Eu circle (8pt);
    \filldraw (\a,\a)++\Eu++\Nd circle (8pt);
    \filldraw (\a,\a)++\Eu++\Nd++\Nu circle (8pt);
    \filldraw (\a,\a)++\Eu++\Nd++\Nu++\Nd circle (8pt);
    \filldraw (\a,\a)++\Eu++\Nd++\Nu++\Nd++\Nu circle (8pt);
    \filldraw (\a,\a)++\Eu++\Nd++\Nu++\Nd++\Nu++\Nd circle (8pt);
     \filldraw (\a,\a)++\Eu++\Nd++\Nu++\Nd++\Nu++\Nd++\Nu circle (8pt);
       \filldraw (\a,\a)++\Eu++\Nd++\Nu++\Nd++\Nu++\Nd++\Nu++\Ed circle (8pt);
\end{tikzpicture}\quad
\begin{tikzpicture}[scale=0.15]
    \draw (\a,\a)-- ++(20,0);
    \draw[dashed,line width=0.1mm] (\a,2.5)-- ++(20,0);
    \draw[dashed,line width=0.1mm] (\a,4.5)-- ++(20,0);
    \draw[dashed,line width=0.1mm] (\a,6.5)-- ++(20,0);
    \draw[dashed,line width=0.1mm] (\a,-1.5)-- ++(20,0);
    \draw[dashed,line width=0.1mm] (\a,-3.5)-- ++(20,0);
    \draw[dashed,line width=0.1mm] (\a,-5.5)-- ++(20,0);
    \draw (\a,-5.5) -- (\a,6.5);
    \draw[dashed,line width=0.1mm] (2.5,-5.5) --(2.5,6.5);
    \draw[dashed,line width=0.1mm] (4.5,-5.5) --(4.5,6.5);
    \draw[dashed,line width=0.1mm] (6.5,-5.5) --(6.5,6.5);
    \draw[dashed,line width=0.1mm] (8.5,-5.5) --(8.5,6.5);
    \draw[dashed,line width=0.1mm] (10.5,-5.5) --(10.5,6.5);
    \draw[dashed,line width=0.1mm] (12.5,-5.5) --(12.5,6.5);
    \draw[dashed,line width=0.1mm] (14.5,-5.5) --(14.5,6.5);
    \draw[dashed,line width=0.1mm] (16.5,-5.5) --(16.5,6.5);
    \draw[dashed,line width=0.1mm] (18.5,-5.5) --(18.5,6.5);
    \draw[dashed,line width=0.1mm] (20.5,-5.5) --(20.5,6.5);
    \draw[solid,line width=0.4mm] (\a,\a)-- ++\Nu-- ++\Ed -- ++\Nu-- ++\Nd -- ++\Nu -- ++\Nd -- ++\Eu -- ++\Nd;
     \filldraw (\a,\a) circle (8pt);
    \filldraw (\a,\a)++\Nu circle (8pt);
    \filldraw (\a,\a)++\Nu++\Ed circle (8pt);
    \filldraw (\a,\a)++\Nu++\Ed++\Nu circle (8pt);
    \filldraw (\a,\a)++\Nu++\Ed++\Nu++\Nd circle (8pt);
    \filldraw (\a,\a)++\Nu++\Ed++\Nu++\Nd++\Nu circle (8pt);
    \filldraw (\a,\a)++\Nu++\Ed++\Nu++\Nd++\Nu++\Nd circle (8pt);
     \filldraw (\a,\a)++\Nu++\Ed++\Nu++\Nd++\Nu++\Nd++\Eu circle (8pt);
       \filldraw (\a,\a)++\Nu++\Ed++\Nu++\Nd++\Nu++\Nd++\Eu++\Nd circle (8pt);
\end{tikzpicture}
\quad
\begin{tikzpicture}[scale=0.15]
    \draw (\a,\a)-- ++(20,0);
    \draw[dashed,line width=0.1mm] (\a,2.5)-- ++(20,0);
    \draw[dashed,line width=0.1mm] (\a,4.5)-- ++(20,0);
    \draw[dashed,line width=0.1mm] (\a,6.5)-- ++(20,0);
    \draw[dashed,line width=0.1mm] (\a,-1.5)-- ++(20,0);
    \draw[dashed,line width=0.1mm] (\a,-3.5)-- ++(20,0);
    \draw[dashed,line width=0.1mm] (\a,-5.5)-- ++(20,0);
    \draw (\a,-5.5) -- (\a,6.5);
    \draw[dashed,line width=0.1mm] (2.5,-5.5) --(2.5,6.5);
    \draw[dashed,line width=0.1mm] (4.5,-5.5) --(4.5,6.5);
    \draw[dashed,line width=0.1mm] (6.5,-5.5) --(6.5,6.5);
    \draw[dashed,line width=0.1mm] (8.5,-5.5) --(8.5,6.5);
    \draw[dashed,line width=0.1mm] (10.5,-5.5) --(10.5,6.5);
    \draw[dashed,line width=0.1mm] (12.5,-5.5) --(12.5,6.5);
    \draw[dashed,line width=0.1mm] (14.5,-5.5) --(14.5,6.5);
    \draw[dashed,line width=0.1mm] (16.5,-5.5) --(16.5,6.5);
    \draw[dashed,line width=0.1mm] (18.5,-5.5) --(18.5,6.5);
    \draw[dashed,line width=0.1mm] (20.5,-5.5) --(20.5,6.5);
    \draw[solid,line width=0.4mm] (\a,\a)-- ++\Nu-- ++\Ed -- ++\Eu-- ++\Ed -- ++\Eu -- ++\Nd;
     \filldraw (\a,\a) circle (8pt);
    \filldraw (\a,\a)++\Nu circle (8pt);
    \filldraw (\a,\a)++\Nu++\Ed circle (8pt);
    \filldraw (\a,\a)++\Nu++\Ed++\Eu circle (8pt);
    \filldraw (\a,\a)++\Nu++\Ed++\Eu++\Ed circle (8pt);
    \filldraw (\a,\a)++\Nu++\Ed++\Eu++\Ed++\Eu circle (8pt);
    \filldraw (\a,\a)++\Nu++\Ed++\Eu++\Ed++\Eu++\Nd circle (8pt);
\end{tikzpicture}
    \caption{The six grand zigzag knight's paths of size 10 having exactly the starting and the ending point on the $x$-axis.}\label{zigzag in B}
    \end{center}
   \end{figure}

The first terms of $b_n$ for $0\leq n\leq 22$ are
$$1,\, 0,\,2, \,0, \,2,\, 2, \,4, \,2, \,4, \,2, \,\mathbf{6}, \,2, \,10, \, 2, \, 18, \, 2, \, 36, \, 2, \, 76, \, 2, \, 166, \, 2, \, 372.$$
    All the odd terms equal 2 for $2n+1\geq 5$, because the only paths of size $2n+1$ in $\mathcal{B}$ are $E\Bar{N}N\cdots N\Bar{N}E$ and its symmetric $\Bar{E}N\Bar{N} \cdots \Bar{N}N\Bar{E}$. Figure~\ref{zigzag in B} shows the six grand zigzag knight paths of size 10 in $\mathcal{B}$ (see the entry indicated by the boldface in the list above).

\subsection{A bijective approach}
In this part, we exhibit a bijection between pairs of integer compositions and grand zigzag knight's paths ending at $(n,k)$ where $n$ and $k$ have the same parity.

In \cite{BK}, B\'ona and Knopfmacher provided the following bijection $\varphi$ between pairs of compositions of $n$ with parts in $\{1,2\}$ that have the same number of parts, and lattice paths starting at $(0,0)$ with size $n$ and steps from the set $\{(1,0),(2,0),(2,1),(1,-1)\}$. Let $X=(x_1,\ldots,x_k)$ and $Y=(y_1,\ldots, y_k)$ be compositions  of $n$ with $k$ parts. The mapping is defined by $\varphi(X,Y)=\varphi(x_1,y_1)\cdots \varphi(x_k,y_k)$ with 
$$\varphi(x_i,y_i)=\begin{cases} (1,0),  & \text{ if } x_i=y_i=1, \\
   (2,0), & \text{ if } x_i=y_i=2, \\
   (2,1), &  \text{ if } x_i=2 \text{ and } y_i=1, \\
   (1,-1), & \text{ if } x_i=1 \text{ and } y_i=2.
\end{cases}$$
Here we construct a similar bijection with grand zigzag knight's paths ending at a point whose coordinates have the same parity.

Let $\mathcal{C}_{n,m}$ be the set of ordered pairs $(X, Y)$ of compositions of $n$ and $m$, respectively, where every part is 1 or 2, such that $X$ and $Y$ have the same number of parts. Note that $\mathcal{C}_{n,m}=\varnothing$ if $n\leq \lceil m/2\rceil$ or $m\leq \lceil n/2\rceil$.

\begin{lemma}\label{cardinal C}
    If $n\leq m$, then $$|\mathcal{C}_{m,n}|=|\mathcal{C}_{n,m}|=\sum_{i=0}^{n-\lceil m/2\rceil}\binom{n-i}{i}\binom{n-i}{m-n+i}.$$
\end{lemma}
\begin{proof}
    Let us count the compositions of $n$ into parts equal to 1 or 2 with $i$ parts ($n\geq i\geq \lceil n/2\rceil$). Since there are $i$ parts, there are necessarily $n-i$ parts equal to 2, and $2i-n$ parts equal to 1. Then it remains to choose the places of the 1's among those $i$ parts, which can be done in $\binom{i}{2i-n}=\binom{i}{n-i}$ ways. Therefore, if we assume $n\leq m$, \[|\mathcal{C}_{n,m}|=\sum_{i=\lceil m/2\rceil}^n \binom{i}{n-i}\binom{i}{m-i}=\sum_{i=0}^{n-\lceil m/2\rceil}\binom{n-i}{i}\binom{n-i}{m-n+i}. \qedhere\]
\end{proof}

\begin{lemma}\label{Bijection C Z}
    If $n,k\in\N$ with $n=k \Mod{2}$, then there is a bijection $\phi$ between $\mathcal{C}_{\frac{n-k}{2},\frac{n+k}{2}}$ and $\mathcal{Z}_{n,k}^+$, and a bijection $\psi$ between $\mathcal{C}_{\frac{n-k}{2},\frac{n+k}{2}}$ and $\mathcal{Z}_{n,k}^-$.
\end{lemma}
\begin{proof}
    Let $(X,Y)\in \mathcal{C}_{\frac{n-k}{2},\frac{n+k}{2}}$. Let $i$ be the number of parts of $X$ and $Y$, so that $X=(x_1,\ldots,x_i)$ and $Y=(y_1,\ldots,y_i)$, with $x_j,y_j\in\{1,2\}$. We define $\phi(X,Y)$ as the path $\phi(x_1,y_1)\cdots \phi(x_k,y_k)$, where
    $$\phi(x_j,y_j)=\begin{cases}
        E\bar{E},   &  \text{ if } x_j=y_j=2,\\
        N\bar{N},   &  \text{ if } x_j=y_j=1,\\
        N\bar{E},   &  \text{ if } x_j=1 \text{ and } y_j=2,\\
        E\bar{N},   &  \text{ if } x_j=2 \text{ and } y_j=1.\\
    \end{cases}$$
    Note that the size of $\phi(x_j,y_j)$ is $x_j+y_j$, so the size of $\phi(X,Y)$ is $\sum_{j=1}^i (x_j+y_j)=n$. The altitude of $\phi(x_j,y_j)$ is $y_j-x_j$, and consequently, the altitude of $\phi(X,Y)$ is $\sum_{j=0}^i (y_j-x_j)=k$. Therefore, $\phi$ maps into $\mathcal{Z}_{n,k}^+$. The inverse of $\phi$ is easy to obtain, since each path of $\mathcal{Z}_{n,k}^+$ can be seen as a path of steps belonging to $\{E\bar{E},N\bar{N},N\bar{E},E\bar{N}\}$. Thus, $\phi$ is a bijection between $\mathcal{C}_{n,m}$ and $\mathcal{Z}_{n,k}^+$. The bijection $\psi$ can be constructed similarly since each path of $\mathcal{Z}_{n,k}^-$ can be seen as a path of steps belonging to $\{\bar{E}E,\bar{N}N,\bar{E}N,\bar{N}E\}$.
\end{proof}

\begin{example}
    If we consider the compositions $$X=(2,2,2,1,1,1,1,2,1) \quad \text{and} \quad Y=(1,2,1,2,2,1,2,1,2),$$
    then $(X,Y)\in\mathcal{C}_{13,14}$ is mapped to the path $E\bar{N}E\bar{E}E\bar{N}N\bar{E}N\bar{E}N\bar{N}N\bar{E}E\bar{N}N\bar{E}\in\mathcal{Z}_{27,1}^+$, see Figure \ref{example bijection}.
    \begin{figure}[ht]
  \begin{center}
         \begin{tikzpicture}[scale=0.15]
    \draw (\a,\a)-- ++(54,0);
    \draw[dashed,line width=0.1mm] (\a,2.5)-- ++(54,0);
    \draw[dashed,line width=0.1mm] (\a,4.5)-- ++(54,0);
    \draw[dashed,line width=0.1mm] (\a,6.5)-- ++(54,0);
    \draw[dashed,line width=0.1mm] (\a,-1.5)-- ++(54,0);
    \draw[dashed,line width=0.1mm] (\a,-3.5)-- ++(54,0);
    \draw[dashed,line width=0.1mm] (\a,-5.5)-- ++(54,0);
    \draw (\a,-5.5) -- (\a,6.5);
    \draw[dashed,line width=0.1mm] (2.5,-5.5) --(2.5,6.5);
    \draw[dashed,line width=0.1mm] (4.5,-5.5) --(4.5,6.5);
    \draw[dashed,line width=0.1mm] (6.5,-5.5) --(6.5,6.5);
    \draw[dashed,line width=0.1mm] (8.5,-5.5) --(8.5,6.5);
    \draw[dashed,line width=0.1mm] (10.5,-5.5) --(10.5,6.5);
    \draw[dashed,line width=0.1mm] (12.5,-5.5) --(12.5,6.5);
    \draw[dashed,line width=0.1mm] (14.5,-5.5) --(14.5,6.5);
    \draw[dashed,line width=0.1mm] (16.5,-5.5) --(16.5,6.5);
    \draw[dashed,line width=0.1mm] (18.5,-5.5) --(18.5,6.5);
    \draw[dashed,line width=0.1mm] (20.5,-5.5) --(20.5,6.5);
    \draw[dashed,line width=0.1mm] (22.5,-5.5) --(22.5,6.5);
    \draw[dashed,line width=0.1mm] (24.5,-5.5) --(24.5,6.5);
    \draw[dashed,line width=0.1mm] (26.5,-5.5) --(26.5,6.5);
    \draw[dashed,line width=0.1mm] (28.5,-5.5) --(28.5,6.5);
    \draw[dashed,line width=0.1mm] (30.5,-5.5) --(30.5,6.5);
    \draw[dashed,line width=0.1mm] (32.5,-5.5) --(32.5,6.5);
    \draw[dashed,line width=0.1mm] (34.5,-5.5) --(34.5,6.5);
    \draw[dashed,line width=0.1mm] (36.5,-5.5) --(36.5,6.5);
    \draw[dashed,line width=0.1mm] (38.5,-5.5) --(38.5,6.5);
    \draw[dashed,line width=0.1mm] (40.5,-5.5) --(40.5,6.5);
    \draw[dashed,line width=0.1mm] (42.5,-5.5) --(42.5,6.5);
    \draw[dashed,line width=0.1mm] (44.5,-5.5) --(44.5,6.5);
    \draw[dashed,line width=0.1mm] (46.5,-5.5) --(46.5,6.5);
    \draw[dashed,line width=0.1mm] (48.5,-5.5) --(48.5,6.5);
    \draw[dashed,line width=0.1mm] (50.5,-5.5) --(50.5,6.5);
    \draw[dashed,line width=0.1mm] (52.5,-5.5) --(52.5,6.5);
    \draw[dashed,line width=0.1mm] (54.5,-5.5) --(54.5,6.5);
    \draw[solid,line width=0.4mm] (\a,\a)-- ++\Eu  -- ++\Nd-- ++\Eu-- ++\Ed--++\Eu-- ++\Nd-- ++\Nu-- ++\Ed-- ++\Nu-- ++\Ed-- ++\Nu-- ++\Nd-- ++\Nu-- ++\Ed-- ++\Eu-- ++\Nd-- ++\Nu-- ++\Ed;
     \filldraw (\a,\a) circle (8pt);
    \filldraw (\a,\a)++\Eu circle (8pt);
    \filldraw (\a,\a)++\Eu++\Nd circle (8pt);
    \filldraw (\a,\a)++\Eu++\Nd++\Eu circle (8pt);
    \filldraw (\a,\a)++\Eu++\Nd++\Eu++\Ed circle (8pt);
    \filldraw (\a,\a)++\Eu++\Nd++\Eu++\Ed++\Eu circle (8pt);
    \filldraw (\a,\a)++\Eu++\Nd++\Eu++\Ed++\Eu++\Nd circle (8pt);
       \filldraw (\a,\a)++\Eu++\Nd++\Eu++\Ed++\Eu++\Nd++\Nu circle (8pt);
        \filldraw (\a,\a)++\Eu++\Nd++\Eu++\Ed++\Eu++\Nd++\Nu++\Ed circle (8pt);
          \filldraw (\a,\a)++\Eu++\Nd++\Eu++\Ed++\Eu++\Nd++\Nu++\Ed++\Nu circle (8pt);
    \filldraw (\a,\a)++\Eu++\Nd++\Eu++\Ed++\Eu++\Nd++\Nu++\Ed++\Nu++\Ed circle (8pt);
    \filldraw (\a,\a)++\Eu++\Nd++\Eu++\Ed++\Eu++\Nd++\Nu++\Ed++\Nu++\Ed++\Nu circle (8pt);
     \filldraw (\a,\a)++\Eu++\Nd++\Eu++\Ed++\Eu++\Nd++\Nu++\Ed++\Nu++\Ed++\Nu++\Nd circle (8pt);
      \filldraw (\a,\a)++\Eu++\Nd++\Eu++\Ed++\Eu++\Nd++\Nu++\Ed++\Nu++\Ed++\Nu++\Nd++\Nu circle (8pt);
       \filldraw (\a,\a)++\Eu++\Nd++\Eu++\Ed++\Eu++\Nd++\Nu++\Ed++\Nu++\Ed++\Nu++\Nd++\Nu++\Ed circle (8pt);
       \filldraw (\a,\a)++\Eu++\Nd++\Eu++\Ed++\Eu++\Nd++\Nu++\Ed++\Nu++\Ed++\Nu++\Nd++\Nu++\Ed++\Eu circle (8pt);
       \filldraw (\a,\a)++\Eu++\Nd++\Eu++\Ed++\Eu++\Nd++\Nu++\Ed++\Nu++\Ed++\Nu++\Nd++\Nu++\Ed++\Eu++\Nd circle (8pt);
       \filldraw (\a,\a)++\Eu++\Nd++\Eu++\Ed++\Eu++\Nd++\Nu++\Ed++\Nu++\Ed++\Nu++\Nd++\Nu++\Ed++\Eu++\Nd++\Nu circle (8pt);
       \filldraw (\a,\a)++\Eu++\Nd++\Eu++\Ed++\Eu++\Nd++\Nu++\Ed++\Nu++\Ed++\Nu++\Nd++\Nu++\Ed++\Eu++\Nd++\Nu++\Ed circle (8pt);
\end{tikzpicture}
\end{center}
\caption{Illustration to Example~3.8. The path $\phi(X,Y)$ is a grand zigzag knight's path of size 27 and altitude~1.}
          \label{example bijection}
\end{figure}
\end{example}

\begin{theorem}\label{closed form zigzag}
    If $n=k \Mod{2}$, with $(n,k)\ne (0,0)$, then $$|\mathcal{Z}_{n,k}|=2\sum_{i=0}^{\frac{n-|k|}{2}}\binom{\frac{n-|k|}{2}-i}{i}\binom{\frac{n-|k|}{2}-i}{|k|+i}.$$

   If $n\ne k \Mod{2}$, then  $$|\mathcal{Z}_{n,k}|=|\mathcal{Z}_{n-1,k-2}^+|+|\mathcal{Z}_{n-2,k-1}^+|+|\mathcal{Z}_{n-1,k+2}^-|+|\mathcal{Z}_{n-2,k+1}^-|,$$ where each of those four terms can be expressed  with the above formula (without the factor 2).
\end{theorem}
\begin{proof}
By Lemma~\ref{Bijection C Z}, if $n=k \Mod{2}$, then we have $|\mathcal{Z}_{n,k}|=|\mathcal{Z}_{n,k}^+|+|\mathcal{Z}_{n,k}^-|=2|\mathcal{C}_{\frac{n-k}{2},\frac{n+k}{2}}|$. We conclude using Lemma~\ref{cardinal C}. If $n\ne k \Mod{2}$, then we have $$\mathcal{Z}_{n,k}=\mathcal{Z}_{n-1,k-2}\cdot N\cup\mathcal{Z}_{n-2,k-1}\cdot E\cup\mathcal{Z}_{n-1,k+2}\cdot\bar{N}\cup\mathcal{Z}_{n-2,k+1}\cdot\bar{E},$$ which completes the proof.
\end{proof}

\begin{rmk}
    When $n=k \Mod{2}$, $|\mathcal{Z}_{n,k}|$ is even (except for $(n,k)=(0,0)$). Indeed, since each path of $\mathcal{Z}_{n,k}$ has an even number of steps, we have a bijection from $\mathcal{Z}_{n,k}^+$ to $\mathcal{Z}_{n,k}^-$, induced by
    $$\begin{array}{ccc}
       N\bar{N}  & \to & \Bar{N}N\\
       E\Bar{E}  & \to & \Bar{E}E\\
       N\Bar{E} & \to & \Bar{E}N\\
       E\Bar{N} & \to & \Bar{N}E.
    \end{array}$$
\end{rmk}

\subsection{Step number of a grand zigzag knight's path}

For $n/2\leq i\leq n$, let $\mathcal{Z}_{n,k}^i$ be the subset of $\mathcal{Z}_{n,k}$ consisting of paths with $i$ steps (see Section~\ref{unrestricted zigzag} for a definition of $\mathcal{Z}_{n,k}$). We define similarly $\mathcal{Z}_{n,k}^{i,+}$ and $\mathcal{Z}_{n,k}^{i,-}$. Note that $|\mathcal{Z}_{n,k}^i|=0$ if $i\ne n-k \Mod{2}$.

\begin{theorem}
    For $n=k \Mod{2}$ and $i$ even, $$|\mathcal{Z}_{n,k}^{i,+}|=|\mathcal{Z}_{n,k}^{i,-}|=\binom{i/2}{\frac{n-i-k}{2}}\binom{i/2}{\frac{n-i+k}{2}}.$$
    For $n\ne k \Mod{2}$ and $i$ odd,
    $$|\mathcal{Z}_{n,k}^{i,+}|=\binom{\frac{i+1}{2}}{\frac{n-k-i}{2}+1}\binom{\frac{i-1}{2}}{\frac{n+k-i}{2}-1},$$
    $$|\mathcal{Z}_{n,k}^{i,-}|=\binom{\frac{i-1}{2}}{\frac{n-k-i}{2}-1}\binom{\frac{i+1}{2}}{\frac{n+k-i}{2}+1}.$$
\end{theorem}
\begin{proof}
    Let us first assume $n=k \Mod{2}$ and $i$ is even. Let $\textbf{e},\textbf{n},\Bar{\textbf{e}},\Bar{\textbf{n}}$ be the numbers of steps $E,N,\Bar{E},\Bar{N}$ of a path in $\mathcal{Z}_{n,k}^{i,+}$ (or $\mathcal{Z}_{n,k}^{i,-}$). Then we have the following equations: 
    $$\left\{
    \begin{array}{ccc}
        \textbf{e} + \Bar{\textbf{e}} & = & n-i , \\
        \textbf{n} + \Bar{\textbf{n}} & = & 2i-n ,\\
        \textbf{e} + \textbf{n} & = & i/2, \\
        \Bar{\textbf{e}} + \Bar{\textbf{n}} & = & i/2, \\
        \textbf{e} + 2\textbf{n} - \Bar{\textbf{e}} - 2\Bar{\textbf{n}} & = & k.
    \end{array}
    \right.$$
    This linear system is non-degenerate, thus we deduce the unique solution 
    $$(\textbf{e},\textbf{n},\Bar{\textbf{e}},\Bar{\textbf{n}})=\left(\frac{n-i-k}{2},\frac{2i-n+k}{2},\frac{n-i+k}{2},\frac{2i-n-k}{2}\right).$$
    Then, it remains to choose the positions of the $N$'s among the up-steps, and the positions of the $\Bar{N}$'s among the down-steps. This eventually yields the desired formula. When $n\ne k \Mod{2}$ and $i$ is odd, the third and the fourth equations become 
    $$\left\{\begin{array}{ccc}
        \textbf{e}+\textbf{n} & = & (i+1)/2,  \\
         \Bar{\textbf{e}} + \Bar{\textbf{n}} & = & (i-1)/2,
    \end{array}\right.$$
     for $\mathcal{Z}_{n,k}^{i,+}$, and 
     $$\left\{\begin{array}{ccc}
        \textbf{e}+\textbf{n} & = & (i-1)/2 , \\
         \Bar{\textbf{e}} + \Bar{\textbf{n}} & = & (i+1)/2,
    \end{array}\right.$$
     for $\mathcal{Z}_{n,k}^{i,-}$. We then conclude similarly.
\end{proof}

\begin{corollary}
    If $n=k \Mod{2}$ with $(n,k)\ne (0,0)$,
    $$|\mathcal{Z}_{n,k}|=2\sum_{\substack{i=0\\i \text{ even}}}^{n-k}\binom{i/2}{\frac{n-i-k}{2}}\binom{i/2}{\frac{n-i+k}{2}}.$$
    If $n\ne k \Mod{2}$,
     $$|\mathcal{Z}_{n,k}|=\sum_{\substack{i=0\\i \text{ odd}}}^{n-k+1}\left[\binom{\frac{i+1}{2}}{\frac{n-k-i}{2}+1}\binom{\frac{i-1}{2}}{\frac{n+k-i}{2}-1}+\binom{\frac{i-1}{2}}{\frac{n-k-i}{2}-1}\binom{\frac{i+1}{2}}{\frac{n+k-i}{2}+1}\right].$$
\end{corollary}

\begin{rmk}
    Note that we obtain the same formula as the one in Theorem~\ref{closed form zigzag}. Also, we now have a formula for $|\mathcal{Z}_{n,k}|$ when $n\ne k \Mod{2}$ involving only 2 sums, instead of 4 as stated in Theorem~\ref{closed form zigzag}.
\end{rmk}

\begin{corollary}\label{expected value nb of steps}
    The expected value for the number of steps of a grand  zigzag knight's path ending at $(n,k)$ is 
    $$\frac{2}{|\mathcal{Z}_{n,k}|}\sum_{\substack{i=0\\i \text{ even}}}^{n-k}i\binom{i/2}{\frac{n-i-k}{2}}\binom{i/2}{\frac{n-i+k}{2}}\quad \mbox{ if }~n=k \Mod{2} ~\mbox{ and}$$
    $$\frac{1}{|\mathcal{Z}_{n,k}|}\sum_{\substack{i=0\\i \text{ odd}}}^{n-k}i\left[\binom{\frac{i+1}{2}}{\frac{n-k-i}{2}+1}\binom{\frac{i-1}{2}}{\frac{n+k-i}{2}-1}+\binom{\frac{i-1}{2}}{\frac{n-k-i}{2}-1}\binom{\frac{i+1}{2}}{\frac{n+k-i}{2}+1}\right],$$
    otherwise.
\end{corollary}

\begin{theorem}
    An asymptotic approximation for the expected number of steps of a grand zigzag knight's path ending on the $x$-axis of size $2n$ (which is also the expected number of parts of a pair of compositions of $n$ with parts in $\{1,2\}$ and having the same number of parts) is  $$\frac{1+\sqrt{5}}{2\sqrt{5}}\cdot 2n.$$
\end{theorem}
\begin{proof}
    By Corollary~\ref{expected value nb of steps} and Theorem~\ref{closed form zigzag}, it suffices to estimate $\sum_{i=0}^n \binom{i}{n-i}^2$ and $\sum_{i=0}^n i\binom{i}{n-i}^2$ (respectively \seqnum{A051286} and \seqnum{A182879}). Their generating functions are respectively  $$\frac{1}{\sqrt{1 - 2z - z^2 - z^3 + z^4}}, \quad \mbox{ and }\quad \frac{z(1+2z^2-z^3)}{((1-3z+z^2)(1+z+z^2))^{3/2}}.$$ We conclude using singularity analysis.
\end{proof}

Based on Mathematica calculations, we suspect that the behavior is the same for odd sizes.
\begin{conjecture}
    An asymptotic approximation for the expected number of steps of a grand zigzag knight's path ending on the $x$-axis of size $n$ is $$\frac{1+\sqrt{5}}{2\sqrt{5}}\cdot n.$$
\end{conjecture}

\section{Bounded grand zigzag knight's paths}\label{bounded grand zigzag}

In this section we focus on grand zigzag knight's paths staying in some regions delimited by horizontal lines.

\subsection{Grand zigzag knight's paths staying above a horizontal line} First, we consider paths that stay above the line $y=-m$ for a given integer $m\geq 0$. By symmetry, this is equivalent to count grand zigzag knight's paths staying below the line $y=+m$. Note that the case $m=0$ has been done in \cite{BarRa}, thus we will assume $m\geq1$. 

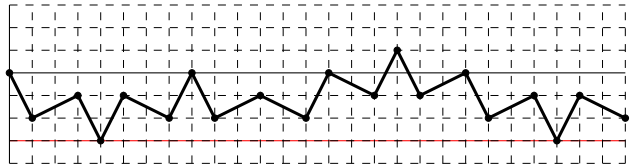
\begin{figure}[ht]
  \begin{center}
         \begin{tikzpicture}[scale=0.15]
    \draw (\a,\a)-- ++(54,0);
    \draw[red] (\a,-5.5)-- ++(54,0);
    \draw[dashed,line width=0.1mm] (\a,2.5)-- ++(54,0);
    \draw[dashed,line width=0.1mm] (\a,4.5)-- ++(54,0);
    \draw[dashed,line width=0.1mm] (\a,6.5)-- ++(54,0);
    \draw[dashed,line width=0.1mm] (\a,-1.5)-- ++(54,0);
    \draw[dashed,line width=0.1mm] (\a,-3.5)-- ++(54,0);
    \draw[dashed,line width=0.1mm] (\a,-5.5)-- ++(54,0);
    \draw[dashed,line width=0.1mm] (\a,-7.5)-- ++(54,0);
    \draw (\a,-7.5) -- (\a,6.5);
    \draw[dashed,line width=0.1mm] (2.5,-7.5) --(2.5,6.5);
    \draw[dashed,line width=0.1mm] (4.5,-7.5) --(4.5,6.5);
    \draw[dashed,line width=0.1mm] (6.5,-7.5) --(6.5,6.5);
    \draw[dashed,line width=0.1mm] (8.5,-7.5) --(8.5,6.5);
    \draw[dashed,line width=0.1mm] (10.5,-7.5) --(10.5,6.5);
    \draw[dashed,line width=0.1mm] (12.5,-7.5) --(12.5,6.5);
    \draw[dashed,line width=0.1mm] (14.5,-7.5) --(14.5,6.5);
    \draw[dashed,line width=0.1mm] (16.5,-7.5) --(16.5,6.5);
    \draw[dashed,line width=0.1mm] (18.5,-7.5) --(18.5,6.5);
    \draw[dashed,line width=0.1mm] (20.5,-7.5) --(20.5,6.5);
    \draw[dashed,line width=0.1mm] (22.5,-7.5) --(22.5,6.5);
    \draw[dashed,line width=0.1mm] (24.5,-7.5) --(24.5,6.5);
    \draw[dashed,line width=0.1mm] (26.5,-7.5) --(26.5,6.5);
    \draw[dashed,line width=0.1mm] (28.5,-7.5) --(28.5,6.5);
    \draw[dashed,line width=0.1mm] (30.5,-7.5) --(30.5,6.5);
    \draw[dashed,line width=0.1mm] (32.5,-7.5) --(32.5,6.5);
    \draw[dashed,line width=0.1mm] (34.5,-7.5) --(34.5,6.5);
    \draw[dashed,line width=0.1mm] (36.5,-7.5) --(36.5,6.5);
    \draw[dashed,line width=0.1mm] (38.5,-7.5) --(38.5,6.5);
    \draw[dashed,line width=0.1mm] (40.5,-7.5) --(40.5,6.5);
    \draw[dashed,line width=0.1mm] (42.5,-7.5) --(42.5,6.5);
    \draw[dashed,line width=0.1mm] (44.5,-7.5) --(44.5,6.5);
    \draw[dashed,line width=0.1mm] (46.5,-7.5) --(46.5,6.5);
    \draw[dashed,line width=0.1mm] (48.5,-7.5) --(48.5,6.5);
    \draw[dashed,line width=0.1mm] (50.5,-7.5) --(50.5,6.5);
    \draw[dashed,line width=0.1mm] (52.5,-7.5) --(52.5,6.5);
    \draw[dashed,line width=0.1mm] (54.5,-7.5) --(54.5,6.5);
    \draw[solid,line width=0.4mm] (\a,\a)-- ++\Nd--++\Eu--++\Nd--++\Nu--++\Ed--++\Nu--++\Nd--++\Eu--++\Ed--++\Nu--++\Ed--++\Nu--++\Nd--++\Eu--++\Nd--++\Eu--++\Nd--++\Nu--++\Ed;
    \filldraw (\a,\a) circle (8pt);
    \filldraw (\a,\a)++\Nd circle (8pt);
    \filldraw (\a,\a)++\Nd++\Eu circle (8pt);
    \filldraw (\a,\a)++\Nd++\Eu++\Nd circle (8pt);
    \filldraw (\a,\a)++\Nd++\Eu++\Nd++\Nu circle (8pt);
    \filldraw (\a,\a)++\Nd++\Eu++\Nd++\Nu++\Ed circle (8pt);
    \filldraw (\a,\a)++\Nd++\Eu++\Nd++\Nu++\Ed++\Nu circle (8pt);
    \filldraw (\a,\a)++\Nd++\Eu++\Nd++\Nu++\Ed++\Nu++\Nd circle (8pt);
    \filldraw (\a,\a)++\Nd++\Eu++\Nd++\Nu++\Ed++\Nu++\Nd++\Eu circle (8pt);
      \filldraw (\a,\a)++\Nd++\Eu++\Nd++\Nu++\Ed++\Nu++\Nd++\Eu++\Ed circle (8pt);
        \filldraw (\a,\a)++\Nd++\Eu++\Nd++\Nu++\Ed++\Nu++\Nd++\Eu++\Ed++\Nu circle (8pt);
         \filldraw (\a,\a)++\Nd++\Eu++\Nd++\Nu++\Ed++\Nu++\Nd++\Eu++\Ed++\Nu++\Ed circle (8pt);
        \filldraw (\a,\a)++\Nd++\Eu++\Nd++\Nu++\Ed++\Nu++\Nd++\Eu++\Ed++\Nu++\Ed++\Nu circle (8pt);
     \filldraw (\a,\a)++\Nd++\Eu++\Nd++\Nu++\Ed++\Nu++\Nd++\Eu++\Ed++\Nu++\Ed++\Nu++\Nd circle (8pt);
      \filldraw (\a,\a)++\Nd++\Eu++\Nd++\Nu++\Ed++\Nu++\Nd++\Eu++\Ed++\Nu++\Ed++\Nu++\Nd++\Eu circle (8pt);
      \filldraw (\a,\a)++\Nd++\Eu++\Nd++\Nu++\Ed++\Nu++\Nd++\Eu++\Ed++\Nu++\Ed++\Nu++\Nd++\Eu++\Nd circle (8pt);
        \filldraw (\a,\a)++\Nd++\Eu++\Nd++\Nu++\Ed++\Nu++\Nd++\Eu++\Ed++\Nu++\Ed++\Nu++\Nd++\Eu++\Nd++\Eu circle (8pt);
            \filldraw (\a,\a)++\Nd++\Eu++\Nd++\Nu++\Ed++\Nu++\Nd++\Eu++\Ed++\Nu++\Ed++\Nu++\Nd++\Eu++\Nd++\Eu++\Nd circle (8pt);
             \filldraw (\a,\a)++\Nd++\Eu++\Nd++\Nu++\Ed++\Nu++\Nd++\Eu++\Ed++\Nu++\Ed++\Nu++\Nd++\Eu++\Nd++\Eu++\Nd++\Nu circle (8pt);
              \filldraw (\a,\a)++\Nd++\Eu++\Nd++\Nu++\Ed++\Nu++\Nd++\Eu++\Ed++\Nu++\Ed++\Nu++\Nd++\Eu++\Nd++\Eu++\Nd++\Nu++\Ed circle (8pt);
      \end{tikzpicture}
\end{center}
\caption{A grand zigzag knight's path staying above the line $y=-3$, with size 27, and altitude -2.}
    \label{zigzag staying above a line}
\end{figure}

\begin{rmk}
    When $n=k \Mod{2}$, there is a bijection between those paths ending at $(n,k)$ and pairs of compositions $(X,Y)\in\mathcal{C}_{\frac{n-k}{2},\frac{n+k}{2}}$ such that for all $j$,
    $$-m\leq\sum_{i=1}^j (y_i-x_i)$$
    (see  the bijections $\phi$ and $\psi$ used in Lemma~\ref{Bijection C Z}).
\end{rmk} 
 Figure~\ref{zigzag staying above a line} shows an example of a grand zigzag knight's path above the line $y=-3$, with size 27, and altitude~$-2$.

As in the previous section, for $k\geq0$, let $f_{-m+k}$ (resp.~$g_{-m+k}$) be the generating function of the number of grand zigzag knight's paths staying above the $y$-coordinate $-m$ and ending at altitude $-m+k$ with an up-step $N$ or $E$ (resp.~with a down-step $\bar{N}$ or $\bar{E}$). Then we easily obtain the following equations
\begin{equation}\label{f_{-m},f_{-m+1}}
f_{-m}=0, \ f_{-m+1}=\mathbb{1}_{[m=1]}+z^2g_{-m},
\end{equation}
\begin{equation}\label{f_{-m+k}}
    f_{-m+k}=zg_{-m+k-2}+z^2g_{-m+k-1}+\mathbb{1}_{[k=m]}+\mathbb{1}_{[k=m+1]}z^2+\mathbb{1}_{[k=m+2]}z \quad \text{for } k\geq2
    \end{equation}
and
\begin{equation}\label{g_{-m+k}}
g_{-m+k}=zf_{-m+k+2}+z^2f_{-m+k+1} \text{ for } k\geq0,
\end{equation}
where $\mathbb{1}_{[a=b]}$ is  $1$ if $a=b$, and $0$ otherwise.

Setting $F(u)=\sum_{k\geq0}f_{-m+k}u^k$ and $G(u)=\sum_{k\geq0}g_{-m+k}u^k$, we obtain
\begin{align*}
    F(u)&=zu^{m+1}(z+u)+u^m+zu(z+u)G(u),\\
    G(u)&=-\frac{z}{u}f_{-m+1}+\left(\frac{z}{u^2}+\frac{z^2}{u}\right)F(u).
\end{align*}

Using the kernel method (the roots $r$ and $s$ of the kernel are the same as  for Section~\ref{algmeth}),  we find $$f_{-m+1}=(1+zr)r^{m-1}+\frac{r^{m+1}}{z}.$$

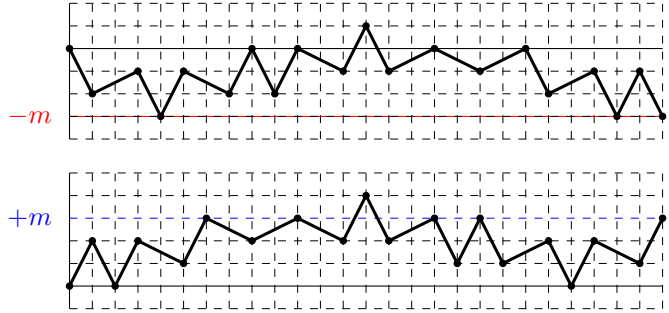
\begin{figure}[ht]
  \begin{center}
         \begin{tikzpicture}[scale=0.15]
    \draw (\a,\a)-- ++(52,0);
    \draw[red] (\a,-5.5)-- ++(52,0);
    \draw[dashed,line width=0.1mm] (\a,2.5)-- ++(52,0);
    \draw[dashed,line width=0.1mm] (\a,4.5)-- ++(52,0);
    \draw[dashed,line width=0.1mm] (\a,-1.5)-- ++(52,0);
    \draw[dashed,line width=0.1mm] (\a,-3.5)-- ++(52,0);
    \draw[dashed,line width=0.1mm] (\a,-5.5)-- ++(52,0);
    \draw[dashed,line width=0.1mm] (\a,-7.5)-- ++(52,0);
    \draw (\a,-7.5) -- (\a,4.5);
    \draw[dashed,line width=0.1mm] (2.5,-7.5) --(2.5,4.5);
    \draw[dashed,line width=0.1mm] (4.5,-7.5) --(4.5,4.5);
    \draw[dashed,line width=0.1mm] (6.5,-7.5) --(6.5,4.5);
    \draw[dashed,line width=0.1mm] (8.5,-7.5) --(8.5,4.5);
    \draw[dashed,line width=0.1mm] (10.5,-7.5) --(10.5,4.5);
    \draw[dashed,line width=0.1mm] (12.5,-7.5) --(12.5,4.5);
    \draw[dashed,line width=0.1mm] (14.5,-7.5) --(14.5,4.5);
    \draw[dashed,line width=0.1mm] (16.5,-7.5) --(16.5,4.5);
    \draw[dashed,line width=0.1mm] (18.5,-7.5) --(18.5,4.5);
    \draw[dashed,line width=0.1mm] (20.5,-7.5) --(20.5,4.5);
    \draw[dashed,line width=0.1mm] (22.5,-7.5) --(22.5,4.5);
    \draw[dashed,line width=0.1mm] (24.5,-7.5) --(24.5,4.5);
    \draw[dashed,line width=0.1mm] (26.5,-7.5) --(26.5,4.5);
    \draw[dashed,line width=0.1mm] (28.5,-7.5) --(28.5,4.5);
    \draw[dashed,line width=0.1mm] (30.5,-7.5) --(30.5,4.5);
    \draw[dashed,line width=0.1mm] (32.5,-7.5) --(32.5,4.5);
    \draw[dashed,line width=0.1mm] (34.5,-7.5) --(34.5,4.5);
    \draw[dashed,line width=0.1mm] (36.5,-7.5) --(36.5,4.5);
    \draw[dashed,line width=0.1mm] (38.5,-7.5) --(38.5,4.5);
    \draw[dashed,line width=0.1mm] (40.5,-7.5) --(40.5,4.5);
    \draw[dashed,line width=0.1mm] (42.5,-7.5) --(42.5,4.5);
    \draw[dashed,line width=0.1mm] (44.5,-7.5) --(44.5,4.5);
    \draw[dashed,line width=0.1mm] (46.5,-7.5) --(46.5,4.5);
    \draw[dashed,line width=0.1mm] (48.5,-7.5) --(48.5,4.5);
    \draw[dashed,line width=0.1mm] (50.5,-7.5) --(50.5,4.5);
    \draw[dashed,line width=0.1mm] (52.5,-7.5) --(52.5,4.5);
    \node at (-3,-5.5) {\color{red} $-m$};
    \draw[solid,line width=0.4mm] (\a,\a)-- ++\Nd--++\Eu--++\Nd--++\Nu--++\Ed--++\Nu--++\Nd--++\Nu--++\Ed--++\Nu--++\Nd--++\Eu--++\Ed--++\Eu--++\Nd--++\Eu--++\Nd--++\Nu--++\Nd;
     \filldraw (\a,\a)++\Nd circle (8pt); \filldraw (\a,\a) circle (8pt);
    \filldraw (\a,\a)++\Nd++\Eu circle (8pt);
    \filldraw (\a,\a)++\Nd++\Eu++\Nd circle (8pt);
    \filldraw (\a,\a)++\Nd++\Eu++\Nd++\Nu circle (8pt);
    \filldraw (\a,\a)++\Nd++\Eu++\Nd++\Nu++\Ed circle (8pt);
     \filldraw (\a,\a)++\Nd++\Eu++\Nd++\Nu++\Ed++\Nu circle (8pt);
     \filldraw (\a,\a)++\Nd++\Eu++\Nd++\Nu++\Ed++\Nu++\Nd circle (8pt);
     \filldraw (\a,\a)++\Nd++\Eu++\Nd++\Nu++\Ed++\Nu++\Nd++\Nu circle (8pt);
     \filldraw (\a,\a)++\Nd++\Eu++\Nd++\Nu++\Ed++\Nu++\Nd++\Nu++\Ed circle (8pt);
      \filldraw (\a,\a)++\Nd++\Eu++\Nd++\Nu++\Ed++\Nu++\Nd++\Nu++\Ed++\Nu circle (8pt);
      \filldraw (\a,\a)++\Nd++\Eu++\Nd++\Nu++\Ed++\Nu++\Nd++\Nu++\Ed++\Nu++\Nd circle (8pt);
     \filldraw (\a,\a)++\Nd++\Eu++\Nd++\Nu++\Ed++\Nu++\Nd++\Nu++\Ed++\Nu++\Nd++\Eu circle (8pt);
      \filldraw (\a,\a)++\Nd++\Eu++\Nd++\Nu++\Ed++\Nu++\Nd++\Nu++\Ed++\Nu++\Nd++\Eu++\Ed circle (8pt);
      \filldraw (\a,\a)++\Nd++\Eu++\Nd++\Nu++\Ed++\Nu++\Nd++\Nu++\Ed++\Nu++\Nd++\Eu++\Ed++\Eu circle (8pt);
      \filldraw (\a,\a)++\Nd++\Eu++\Nd++\Nu++\Ed++\Nu++\Nd++\Nu++\Ed++\Nu++\Nd++\Eu++\Ed++\Eu++\Nd circle (8pt);
       \filldraw (\a,\a)++\Nd++\Eu++\Nd++\Nu++\Ed++\Nu++\Nd++\Nu++\Ed++\Nu++\Nd++\Eu++\Ed++\Eu++\Nd++\Eu circle (8pt);
         \filldraw (\a,\a)++\Nd++\Eu++\Nd++\Nu++\Ed++\Nu++\Nd++\Nu++\Ed++\Nu++\Nd++\Eu++\Ed++\Eu++\Nd++\Eu++\Nd circle (8pt);
         \filldraw (\a,\a)++\Nd++\Eu++\Nd++\Nu++\Ed++\Nu++\Nd++\Nu++\Ed++\Nu++\Nd++\Eu++\Ed++\Eu++\Nd++\Eu++\Nd++\Nu circle (8pt);
          \filldraw (\a,\a)++\Nd++\Eu++\Nd++\Nu++\Ed++\Nu++\Nd++\Nu++\Ed++\Nu++\Nd++\Eu++\Ed++\Eu++\Nd++\Eu++\Nd++\Nu++\Nd circle (8pt);
\end{tikzpicture}\\\vskip0.4cm
\begin{tikzpicture}[scale=0.15]
    \draw (\a,\a)-- ++(52,0);
    \draw[dashed,line width=0.1mm] (\a,2.5)-- ++(52,0);
    \draw[dashed,line width=0.1mm] (\a,4.5)-- ++(52,0);
    \draw[blue,dashed,line width=0.1mm] (\a,6.5)-- ++(52,0);
    \draw[dashed,line width=0.1mm] (\a,8.5)-- ++(52,0);
    \draw[dashed,line width=0.1mm] (\a,10.5)-- ++(52,0);
    \draw[dashed,line width=0.1mm] (\a,-1.5)-- ++(52,0);
    \draw (\a,-1.5) -- (\a,10.5);
    \draw[dashed,line width=0.1mm] (2.5,-1.5) --(2.5,10.5);
    \draw[dashed,line width=0.1mm] (4.5,-1.5) --(4.5,10.5);
    \draw[dashed,line width=0.1mm] (6.5,-1.5) --(6.5,10.5);
    \draw[dashed,line width=0.1mm] (8.5,-1.5) --(8.5,10.5);
    \draw[dashed,line width=0.1mm] (10.5,-1.5) --(10.5,10.5);
    \draw[dashed,line width=0.1mm] (12.5,-1.5) --(12.5,10.5);
    \draw[dashed,line width=0.1mm] (14.5,-1.5) --(14.5,10.5);
    \draw[dashed,line width=0.1mm] (16.5,-1.5) --(16.5,10.5);
    \draw[dashed,line width=0.1mm] (18.5,-1.5) --(18.5,10.5);
    \draw[dashed,line width=0.1mm] (20.5,-1.5) --(20.5,10.5);
    \draw[dashed,line width=0.1mm] (22.5,-1.5) --(22.5,10.5);
    \draw[dashed,line width=0.1mm] (24.5,-1.5) --(24.5,10.5);
    \draw[dashed,line width=0.1mm] (26.5,-1.5) --(26.5,10.5);
    \draw[dashed,line width=0.1mm] (28.5,-1.5) --(28.5,10.5);
    \draw[dashed,line width=0.1mm] (30.5,-1.5) --(30.5,10.5);
    \draw[dashed,line width=0.1mm] (32.5,-1.5) --(32.5,10.5);
    \draw[dashed,line width=0.1mm] (34.5,-1.5) --(34.5,10.5);
    \draw[dashed,line width=0.1mm] (36.5,-1.5) --(36.5,10.5);
    \draw[dashed,line width=0.1mm] (38.5,-1.5) --(38.5,10.5);
    \draw[dashed,line width=0.1mm] (40.5,-1.5) --(40.5,10.5);
    \draw[dashed,line width=0.1mm] (42.5,-1.5) --(42.5,10.5);
    \draw[dashed,line width=0.1mm] (44.5,-1.5) --(44.5,10.5);
    \draw[dashed,line width=0.1mm] (46.5,-1.5) --(46.5,10.5);
    \draw[dashed,line width=0.1mm] (48.5,-1.5) --(48.5,10.5);
    \draw[dashed,line width=0.1mm] (50.5,-1.5) --(50.5,10.5);
    \draw[dashed,line width=0.1mm] (52.5,-1.5) --(52.5,10.5);
    \node at (-3,6.5) {\color{blue} $+m$};
    \draw[solid,line width=0.4mm] (\a,\a)-- ++\Nu-- ++\Nd-- ++\Nu-- ++\Ed-- ++\Nu-- ++\Ed-- ++\Eu-- ++\Ed-- ++\Nu-- ++\Nd-- ++\Eu-- ++\Nd-- ++\Nu-- ++\Nd-- ++\Eu-- ++\Nd-- ++\Nu-- ++\Ed-- ++\Nu;
    \filldraw (\a,\a) circle (8pt);
    \filldraw (\a,\a)++\Nu circle (8pt);
    \filldraw (\a,\a)++\Nu++\Nd circle (8pt);
    \filldraw (\a,\a)++\Nu++\Nd++\Nu circle (8pt);
    \filldraw (\a,\a)++\Nu++\Nd++\Nu++\Ed circle (8pt);
    \filldraw (\a,\a)++\Nu++\Nd++\Nu++\Ed++\Nu circle (8pt);
     \filldraw (\a,\a)++\Nu++\Nd++\Nu++\Ed++\Nu++\Ed circle (8pt);
     \filldraw (\a,\a)++\Nu++\Nd++\Nu++\Ed++\Nu++\Ed+\Eu circle (8pt);
     \filldraw (\a,\a)++\Nu++\Nd++\Nu++\Ed++\Nu++\Ed++\Eu++\Ed circle (8pt);
     \filldraw (\a,\a)++\Nu++\Nd++\Nu++\Ed++\Nu++\Ed++\Eu++\Ed++\Nu circle (8pt);
      \filldraw (\a,\a)++\Nu++\Nd++\Nu++\Ed++\Nu++\Ed++\Eu++\Ed++\Nu++\Nd circle (8pt);
      \filldraw (\a,\a)++\Nu++\Nd++\Nu++\Ed++\Nu++\Ed++\Eu++\Ed++\Nu++\Nd++\Eu circle (8pt);
     \filldraw (\a,\a)++\Nu++\Nd++\Nu++\Ed++\Nu++\Ed++\Eu++\Ed++\Nu++\Nd++\Eu++\Nd circle (8pt);
      \filldraw (\a,\a)++\Nu++\Nd++\Nu++\Ed++\Nu++\Ed++\Eu++\Ed++\Nu++\Nd++\Eu++\Nd++\Nu circle (8pt);
      \filldraw (\a,\a)++\Nu++\Nd++\Nu++\Ed++\Nu++\Ed++\Eu++\Ed++\Nu++\Nd++\Eu++\Nd++\Nu++\Nd circle (8pt);
      \filldraw (\a,\a)++\Nu++\Nd++\Nu++\Ed++\Nu++\Ed++\Eu++\Ed++\Nu++\Nd++\Eu++\Nd++\Nu++\Nd++\Eu circle (8pt);
     \filldraw (\a,\a)++\Nu++\Nd++\Nu++\Ed++\Nu++\Ed++\Eu++\Ed++\Nu++\Nd++\Eu++\Nd++\Nu++\Nd++\Eu++\Nd circle (8pt);
     \filldraw (\a,\a)++\Nu++\Nd++\Nu++\Ed++\Nu++\Ed++\Eu++\Ed++\Nu++\Nd++\Eu++\Nd++\Nu++\Nd++\Eu++\Nd++\Nu circle (8pt);
     \filldraw (\a,\a)++\Nu++\Nd++\Nu++\Ed++\Nu++\Ed++\Eu++\Ed++\Nu++\Nd++\Eu++\Nd++\Nu++\Nd++\Eu++\Nd++\Nu++\Ed circle (8pt);
     \filldraw (\a,\a)++\Nu++\Nd++\Nu++\Ed++\Nu++\Ed++\Eu++\Ed++\Nu++\Nd++\Eu++\Nd++\Nu++\Nd++\Eu++\Nd++\Nu++\Ed++\Nu circle (8pt);
\end{tikzpicture}
\end{center}
\caption{There is a bijection between grand zigzag knight's paths staying above $y=-m$ and ending at altitude~$-m$, and grand zigzag knight's paths staying above $y=0$ and ending at altitude $k$.}
    \label{bijection zigzag above a line}
\end{figure}

\begin{rmk}
    From equation~(\ref{f_{-m},f_{-m+1}}), we obtain $g_{-m}=\frac{1+zr}{z^2}r^{m-1}+\frac{r^{m+1}}{z^3}-\mathbb{1}_{[m=1]}\frac{1}{z^2}$. This is the same expression as the generating function for the total number of grand zigzag knight's paths of altitude $m$ and staying above the $x$-axis, see \cite[Theorem 1]{BarRa}. Indeed, reading the paths backwards, we get a bijection between those ending at altitude $m$ and staying above $y=0$, and those ending at altitude  $-m$ and staying above $y=-m$, see Figure~\ref{bijection zigzag above a line}.
\end{rmk}

\begin{theorem}\label{zigzag above a line}
The bivariate generating functions for the number of grand zigzag knight's paths staying above $y=-m$ and ending with respectively an up-step and a down-step with respect to the size and the altitude are

$$F(u)=-\frac{u(u^m(1+z^2u+zu^2)-r^{m-1}z(u+z)(z+rz^2+r^2))}{z^3(u-r)(u-s)},$$    

$$G(u)=\frac{r^{m-1}(z+r^2+z^2r)-zu^{m-1}(1 + zu)(1 + z^2u + zu^2)}{z^3(u-r)(u-s)}.$$
Setting $H(u):=F(u)+G(u)$, the generating function for the total number of grand zigzag knight's paths staying above $y=-m$ with respect to the size is 
$$H(1)=\frac{zr^{m-1}+z^2r^m+r^{m+1}-z^2-z-1}{z^2+z-1}.$$
\end{theorem}

Here are the first terms of $H(1)$ when $m=2$ and $0\leq n\leq 15$:
$$1,\,2,\,4,\,6,\,9,\,15,\,23,\,38,\,58,\,95,\,147,\,239,\,373,\,603,\,947,\,1525.$$

\begin{rmk}\label{convergence line}
    $H(1)$ converges as $m\to\infty$ to $(1+z+z^2)/(1-z-z^2)$, the generating function for the total number of grand zigzag knight's paths (see the beginning of Section~\ref{unrestricted zigzag}). Moreover, the valuation of $r$ is 3, thus we can estimate the rate of convergence:
    $$\text{val}\left(H(1)-\frac{1+z+z^2}{1-z-z^2}\right)= 3m-2.$$
    In particular, every grand zigzag knight's path of size $\leq 3m-3$ stays above $y=-m$, and we know that the path $(\Bar{N}E)^{m-1}\Bar{N}$ has size $3m-2$ and goes below $y=-m$.
\end{rmk}

\begin{corollary}
    For $m\geq 0$, the probability that a grand zigzag knight's path chosen uniformly at random among all grand zigzag knight's paths of size $n$ stays above the line $y=-m$ is asymptotically $c_m/\sqrt{n}$, with 
    $$c_m=\begin{cases}
       \frac{2+\sqrt{5}}{2}\sqrt{\frac{7\sqrt{5}-15}{\pi}}, & \text{ if } m=0,\\
       \frac{4m+3-\sqrt{5}}{4(\sqrt{5}-2)}\sqrt{\frac{7\sqrt{5}-15}{\pi}},  & \text{ if } m\geq1.
    \end{cases}$$
\end{corollary}
\begin{proof}
    Using basic singularity analysis, since the generating function for the total number of grand zigzag knight's paths with respect to the size is $(1+z+z^2)/(1-z-z^2)$, we deduce that the asymptotic for the total number $a(n)$ of grand zigzag knight's paths of size $n$ is $\frac{4}{\sqrt{5}(\sqrt{5}-1)}\alpha^{-n}$, with $\alpha=\frac{\sqrt{5}-1}{2}$. Now we study $H(1)$ (see Theorem~\ref{zigzag above a line}). It turns out that the main singularity of this function (the one with the smallest absolute value) is also $\alpha$. Using the following approximation near $\alpha$:
    $$\frac{r^{j}}{z-\alpha}=\frac{1}{z-\alpha}+\frac{j\sqrt{10-4\sqrt{5}}}{\sqrt{5}-2}(z-\alpha)^{-1/2}+O(1),$$
    we deduce
    $$H(1)\sim\frac{(2m+1-\alpha)\sqrt{10-4\sqrt{5}}}{\sqrt{5}(\sqrt{5}-2)}(z-\alpha)^{-1/2}$$
    as $z\to\alpha$. We thus have the following approximation for the number $a_m(n)$ of grand zigzag knight's paths of size $n$ and staying above $y=-m$ (see for instance \cite[Chapter VI]{analytic combinatorics}) as $n\to\infty$, when $m\geq 1$:
    $$a_m(n)\sim \frac{(2m+1-\alpha)\sqrt{10-4\sqrt{5}}}{\sqrt{5\alpha}(\sqrt{5}-2)}\frac{\alpha^{-n}}{\sqrt{\pi n}}.$$
    By taking the quotient, we obtain the asymptotic for the desired probability as $n\to\infty$:
    $$\frac{a_m(n)}{a(n)}\sim \frac{4m+3-\sqrt{5}}{4(\sqrt{5}-2)}\sqrt{\frac{7\sqrt{5}-15}{\pi n}}.$$
    When $m=0$, we use similarly the generating function from \cite[Corollary 2]{BarRa}.
\end{proof}


\begin{rmk}
    We can deduce from the previous proof (going one order further) that the asymptotic probability that a grand zigzag knight's path of size $n$ chosen uniformly at random has minimal height $-m$ is 
    $$\frac{a_m(n)-a_{m-1}(n)}{a(n)}\sim \begin{cases}
      \frac{2+\sqrt{5}}{2}\sqrt{\frac{7\sqrt{5}-15}{\pi n}},   & \text{ if } m=0,\\
        \frac{5+3\sqrt{5}}{4}\sqrt{\frac{7\sqrt{5}-15}{\pi n}}, & \text{ if } m=1, \\
       (2+\sqrt{5})\sqrt{\frac{7\sqrt{5}-15}{\pi n}},  & \text{ if } m\geq2. 
    \end{cases}$$
    In particular, asymptotically, every minimal height $\leq -2$ has the same probability. By symmetry, this is exactly the same for maximal heights.
\end{rmk}

\subsection{Grand zigzag knight's paths staying in a symmetric tube}
\begin{figure}[ht]
  \begin{center}
         \begin{tikzpicture}[scale=0.15]
    \draw (\a,\a)-- ++(54,0);
    \draw[red] (\a,-3.5)-- ++(54,0);
    \draw[red] (\a,4.5)-- ++(54,0);
    \draw[dashed,line width=0.1mm] (\a,2.5)-- ++(54,0);
    \draw[dashed,line width=0.1mm] (\a,4.5)-- ++(54,0);
    \draw[dashed,line width=0.1mm] (\a,6.5)-- ++(54,0);
    \draw[dashed,line width=0.1mm] (\a,-1.5)-- ++(54,0);
    \draw[dashed,line width=0.1mm] (\a,-3.5)-- ++(54,0);
    \draw[dashed,line width=0.1mm] (\a,-5.5)-- ++(54,0);
    \draw (\a,-5.5) -- (\a,6.5);
    \draw[dashed,line width=0.1mm] (2.5,-5.5) --(2.5,6.5);
    \draw[dashed,line width=0.1mm] (4.5,-5.5) --(4.5,6.5);
    \draw[dashed,line width=0.1mm] (6.5,-5.5) --(6.5,6.5);
    \draw[dashed,line width=0.1mm] (8.5,-5.5) --(8.5,6.5);
    \draw[dashed,line width=0.1mm] (10.5,-5.5) --(10.5,6.5);
    \draw[dashed,line width=0.1mm] (12.5,-5.5) --(12.5,6.5);
    \draw[dashed,line width=0.1mm] (14.5,-5.5) --(14.5,6.5);
    \draw[dashed,line width=0.1mm] (16.5,-5.5) --(16.5,6.5);
    \draw[dashed,line width=0.1mm] (18.5,-5.5) --(18.5,6.5);
    \draw[dashed,line width=0.1mm] (20.5,-5.5) --(20.5,6.5);
    \draw[dashed,line width=0.1mm] (22.5,-5.5) --(22.5,6.5);
    \draw[dashed,line width=0.1mm] (24.5,-5.5) --(24.5,6.5);
    \draw[dashed,line width=0.1mm] (26.5,-5.5) --(26.5,6.5);
    \draw[dashed,line width=0.1mm] (28.5,-5.5) --(28.5,6.5);
    \draw[dashed,line width=0.1mm] (30.5,-5.5) --(30.5,6.5);
    \draw[dashed,line width=0.1mm] (32.5,-5.5) --(32.5,6.5);
    \draw[dashed,line width=0.1mm] (34.5,-5.5) --(34.5,6.5);
    \draw[dashed,line width=0.1mm] (36.5,-5.5) --(36.5,6.5);
    \draw[dashed,line width=0.1mm] (38.5,-5.5) --(38.5,6.5);
    \draw[dashed,line width=0.1mm] (40.5,-5.5) --(40.5,6.5);
    \draw[dashed,line width=0.1mm] (42.5,-5.5) --(42.5,6.5);
    \draw[dashed,line width=0.1mm] (44.5,-5.5) --(44.5,6.5);
    \draw[dashed,line width=0.1mm] (46.5,-5.5) --(46.5,6.5);
    \draw[dashed,line width=0.1mm] (48.5,-5.5) --(48.5,6.5);
    \draw[dashed,line width=0.1mm] (50.5,-5.5) --(50.5,6.5);
    \draw[dashed,line width=0.1mm] (52.5,-5.5) --(52.5,6.5);
    \draw[dashed,line width=0.1mm] (54.5,-5.5) --(54.5,6.5);
    \draw[solid,line width=0.4mm] (\a,\a)-- ++\Nu--++\Ed--++\Eu--++\Nd--++\Eu--++\Nd--++\Eu--++\Nd--++\Nu--++\Nd--++\Nu--++\Nd--++\Eu--++\Ed--++\Nu--++\Ed--++\Nu--++\Nd--++\Eu;
 \filldraw (\a,\a) circle (8pt);
    \filldraw (\a,\a)++\Nu circle (8pt);
    \filldraw (\a,\a)++\Nu++\Ed circle (8pt);
    \filldraw (\a,\a)++\Nu++\Ed++\Eu circle (8pt);
    \filldraw (\a,\a)++\Nu++\Ed++\Eu++\Nd circle (8pt);
    \filldraw (\a,\a)++\Nu++\Ed++\Eu++\Nd++\Eu circle (8pt);
     \filldraw (\a,\a)++\Nu++\Ed++\Eu++\Nd++\Eu++\Nd circle (8pt);
     \filldraw (\a,\a)++\Nu++\Ed++\Eu++\Nd++\Eu++\Nd+\Eu circle (8pt);
     \filldraw (\a,\a)++\Nu++\Ed++\Eu++\Nd++\Eu++\Nd++\Eu++\Nd circle (8pt);
     \filldraw (\a,\a)++\Nu++\Ed++\Eu++\Nd++\Eu++\Nd++\Eu++\Nd++\Nu circle (8pt);
      \filldraw (\a,\a)++\Nu++\Ed++\Eu++\Nd++\Eu++\Nd++\Eu++\Nd++\Nu++\Nd circle (8pt);
      \filldraw (\a,\a)++\Nu++\Ed++\Eu++\Nd++\Eu++\Nd++\Eu++\Nd++\Nu++\Nd++\Nu circle (8pt);
     \filldraw (\a,\a)++\Nu++\Ed++\Eu++\Nd++\Eu++\Nd++\Eu++\Nd++\Nu++\Nd++\Nu++\Nd circle (8pt);
      \filldraw (\a,\a)++\Nu++\Ed++\Eu++\Nd++\Eu++\Nd++\Eu++\Nd++\Nu++\Nd++\Nu++\Nd++\Eu circle (8pt);
      \filldraw (\a,\a)++\Nu++\Ed++\Eu++\Nd++\Eu++\Nd++\Eu++\Nd++\Nu++\Nd++\Nu++\Nd++\Eu++\Ed circle (8pt);
      \filldraw (\a,\a)++\Nu++\Ed++\Eu++\Nd++\Eu++\Nd++\Eu++\Nd++\Nu++\Nd++\Nu++\Nd++\Eu++\Ed++\Nu circle (8pt);
       \filldraw (\a,\a)++\Nu++\Ed++\Eu++\Nd++\Eu++\Nd++\Eu++\Nd++\Nu++\Nd++\Nu++\Nd++\Eu++\Ed++\Nu ++\Ed circle (8pt);
        \filldraw (\a,\a)++\Nu++\Ed++\Eu++\Nd++\Eu++\Nd++\Eu++\Nd++\Nu++\Nd++\Nu++\Nd++\Eu++\Ed++\Nu ++\Ed++\Nu circle (8pt);
         \filldraw (\a,\a)++\Nu++\Ed++\Eu++\Nd++\Eu++\Nd++\Eu++\Nd++\Nu++\Nd++\Nu++\Nd++\Eu++\Ed++\Nu ++\Ed++\Nu++\Nd circle (8pt);
           \filldraw (\a,\a)++\Nu++\Ed++\Eu++\Nd++\Eu++\Nd++\Eu++\Nd++\Nu++\Nd++\Nu++\Nd++\Eu++\Ed++\Nu ++\Ed++\Nu++\Nd++\Eu circle (8pt);
\end{tikzpicture}
\end{center}
\caption{A grand zigzag knight's path staying between the lines $y=-2$ and $y=+2$, with size 27 and altitude~0.}
    \label{zigzag staying above a lineb}
\end{figure}
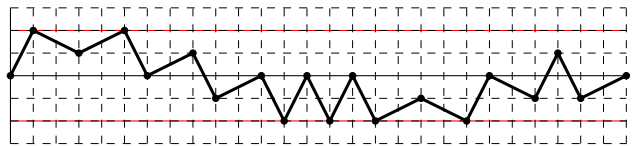

Here we count grand zigzag knight's paths that stay in a symmetric tube, i.e. between two  horizontal lines $y=-m$ and $y=+m$, for $m\geq1$. See~\cite{GuPro,CHNS,Corridor} for previous studies of lattices paths staying in a tube.   Figure~\ref{zigzag staying above a lineb} shows an example of a grand zigzag knight's path staying between the lines $y=-2$ and $y=2$, with size 27 and altitude $0$.

We use the same notation as in the previous subsection. Note that $f_{-k}=g_k$ for all $-m\leq k\leq m$. We obtain the following equations:
\begin{equation}\label{f_{-m},tube sym}
f_{-m}=g_m=0, \ f_{-m+1}= \mathbb{1}_{[m=1]} + z^2g_{-m}, \  g_{m-1}=z^2f_m,
\end{equation}
\begin{equation}\label{f_{-m+k},tube sym}
    f_{-m+k}=zg_{-m+k-2}+z^2g_{-m+k-1}+\mathbb{1}_{[k=m]}+\mathbb{1}_{[k=m+1]}z^2+\mathbb{1}_{[k=m+2]}z
\end{equation}
for $2\leq k \leq 2m$, and
\begin{equation}\label{g_{-m+k},tube sym}
g_{-m+k}=zf_{-m+k+2}+z^2f_{-m+k+1} \text{ \ for } 0\leq k\leq 2m-2.
\end{equation}

\begin{rmk}\label{rational solution}
    We actually have a finite linear system with unknowns $f_k,g_k$, with polynomial coefficients. Therefore, we expect rational solutions and in particular, the generating function for the total number of grand zigzag knight's paths staying between $y=-m$ and $y=+m$ is rational.
\end{rmk}

Setting $F(u)=\sum_{k=0}^{2m}f_{-m+k}u^k$ and $G(u)=\sum_{k=0}^{2m}g_{-m+k}u^k$ we obtain

\begin{align*}
    F(u)&=u^m-zu^{2m+1}(f_{-m+1}-\mathbb{1}_{[m=1]}) + zu(z+u)(G(u)+u^m),\\
    G(u)&=-\frac{z}{u}f_{-m+1}+\left(\frac{z}{u^2}+\frac{z^2}{u}\right)F(u).
\end{align*}
The kernel method then leads to 
$$f_{-m+1}= \frac{r^m(1 + rz^2 + r^2z( 1 +\mathbb{1}_{[m=1]}))}{z(rz + z^2 + r^{2m + 1})}.$$

\begin{theorem}
The generating functions $F(u)$ and $G(u)$ counting the number of grand zigzag knight's paths staying between $y=-m$ and $y=+m$ with respect to the size and the altitude are given by:
\begin{align*}
F(u)&=\scalemath{1}{-\frac{(z(\mathbb{1}_{[m=1]} - f_{-m+1})u^{2m + 1} - z^2(u + z)f_{-m+1} + u^m(1 + z^2u + zu^2))u}{z^3(u-r)(u-s)}},\\
G(u)&=\scalemath{1}{-\frac{z(1 + zu)(\mathbb{1}_{[m=1]} - f_{-m+1})u^{2m + 1} - uf_{-m+1} + u^m(1 + zu)(1 + z^2u + zu^2)}{z^2(u-r)(u-s)}}.
\end{align*}
Setting $H(u)=F(u)+G(u)$, the generating function for the total number of grand zigzag knight's paths staying between $y=-m$ and $y=+m$ is 
$$H(1)=\frac{(2f_{-m+1} -\mathbb{1}_{[m=1]})z - z^2 - z -1}{z^2 + z - 1}.$$
\end{theorem}
\begin{rmk}\label{convergence tube} 
    As a consequence of Remark \ref{rational solution}, $f_{-m+1}$ and $H(1)$ are rational.
    As in Remark~\ref{convergence line}, $H(1)$ converges as $m\to\infty$ to $(1+z+z^2)/(1-z-z^2)$, the generating function for the total number of grand zigzag knight's paths (see (\ref{zsec3}) at the beginning of Section~\ref{unrestricted zigzag}).
\end{rmk}

\begin{example}\label{zigzag small tube}
    For $m=1$, the generating function for the number of grand zigzag knight's paths staying between $y=-1$ and $y=+1$ and ending on the $x$-axis is $(-z^4+z-1)/(z^4+z-1)$. Its first terms for  $0\leq n\leq 18$ are:
    $$1,\, 0,\, 0,\, 0,\, 2,\, 2,\, 2,\, 2,\, 4,\, 6,\, 8,\, 10,\, 14,\, 20,\, 28,\, 38,\, 52,\, 72,\, 100.$$
    The term of order $2n+4$ corresponds to 2$\times$\href{https://oeis.org/A052535}{A052535}($n$)$=2\sum_{k=0}^{\lfloor n/2\rfloor}\binom{2n-3k}{k}$.
    The term of order $2n+3$ corresponds to 2$\times$\href{https://oeis.org/A158943}{A158943}($n$).
\end{example} 
    

    For instance, we can easily exhibit three bijections between such paths from $(0,0)$ to $(2n+4,0)$ starting with $E$ and 
    \begin{itemize}
        \item tilings of a $2\times 2n$ rectangle with $1\times 2$ and $4\times 1$ tiles,
        \item pairs $X=(x_1,\ldots,x_k),Y=(y_1,\ldots,y_k)$ of compositions of $n+2$ with same number of parts, with all parts in $\{1,2\}$ and such that for all $1\leq j\leq k$,
        $$\left|\sum_{i=1}^j (x_i-y_i)\right|\leq 1.$$
        \item compositions of $n$ with parts in $\{2,1,3,5,9,11,\ldots\}$. Indeed, a grand zigzag knight's path of size $2n+4$ staying between $y=-1$ and $y=+1$ can be uniquely decomposed as $E\textbf{P}\Bar{E}$ with $\textbf{P}$ of size $2n$ and having its steps in $\{\Bar{E}E\}\cup\bigcup_{k\geq0}\{\Bar{N}(E\Bar{E})^{k}N\}$. If $\textbf{P}=S_1\cdots S_j$, then we set $\Phi(E\textbf{P}\Bar{E})=(\Phi(S_1),\ldots,\Phi(S_j))$ with
        $\Phi(\Bar{E}E) =  2$ and $\Phi(\Bar{N}(E\Bar{E})^{k}N) = 2k+1$.
    \end{itemize}

\appendix 
\section{Grand zigzag knight's paths staying in a general tube}
\begin{figure}[ht]
  \begin{center}
         \begin{tikzpicture}[scale=0.15]
    \draw (\a,\a)-- ++(54,0);
    \draw[red] (\a,-3.5)-- ++(54,0);
    \draw[red] (\a,6.5)-- ++(54,0);
    \draw[dashed,line width=0.1mm] (\a,2.5)-- ++(54,0);
    \draw[dashed,line width=0.1mm] (\a,4.5)-- ++(54,0);
    \draw[dashed,line width=0.1mm] (\a,6.5)-- ++(54,0);
    \draw[dashed,line width=0.1mm] (\a,8.5)-- ++(54,0);
    \draw[dashed,line width=0.1mm] (\a,-1.5)-- ++(54,0);
    \draw[dashed,line width=0.1mm] (\a,-3.5)-- ++(54,0);
    \draw[dashed,line width=0.1mm] (\a,-5.5)-- ++(54,0);
    \draw (\a,-5.5) -- (\a,8.5);
    \draw[dashed,line width=0.1mm] (2.5,-5.5) --(2.5,8.5);
    \draw[dashed,line width=0.1mm] (4.5,-5.5) --(4.5,8.5);
    \draw[dashed,line width=0.1mm] (6.5,-5.5) --(6.5,8.5);
    \draw[dashed,line width=0.1mm] (8.5,-5.5) --(8.5,8.5);
    \draw[dashed,line width=0.1mm] (10.5,-5.5) --(10.5,8.5);
    \draw[dashed,line width=0.1mm] (12.5,-5.5) --(12.5,8.5);
    \draw[dashed,line width=0.1mm] (14.5,-5.5) --(14.5,8.5);
    \draw[dashed,line width=0.1mm] (16.5,-5.5) --(16.5,8.5);
    \draw[dashed,line width=0.1mm] (18.5,-5.5) --(18.5,8.5);
    \draw[dashed,line width=0.1mm] (20.5,-5.5) --(20.5,8.5);
    \draw[dashed,line width=0.1mm] (22.5,-5.5) --(22.5,8.5);
    \draw[dashed,line width=0.1mm] (24.5,-5.5) --(24.5,8.5);
    \draw[dashed,line width=0.1mm] (26.5,-5.5) --(26.5,8.5);
    \draw[dashed,line width=0.1mm] (28.5,-5.5) --(28.5,8.5);
    \draw[dashed,line width=0.1mm] (30.5,-5.5) --(30.5,8.5);
    \draw[dashed,line width=0.1mm] (32.5,-5.5) --(32.5,8.5);
    \draw[dashed,line width=0.1mm] (34.5,-5.5) --(34.5,8.5);
    \draw[dashed,line width=0.1mm] (36.5,-5.5) --(36.5,8.5);
    \draw[dashed,line width=0.1mm] (38.5,-5.5) --(38.5,8.5);
    \draw[dashed,line width=0.1mm] (40.5,-5.5) --(40.5,8.5);
    \draw[dashed,line width=0.1mm] (42.5,-5.5) --(42.5,8.5);
    \draw[dashed,line width=0.1mm] (44.5,-5.5) --(44.5,8.5);
    \draw[dashed,line width=0.1mm] (46.5,-5.5) --(46.5,8.5);
    \draw[dashed,line width=0.1mm] (48.5,-5.5) --(48.5,8.5);
    \draw[dashed,line width=0.1mm] (50.5,-5.5) --(50.5,8.5);
    \draw[dashed,line width=0.1mm] (52.5,-5.5) --(52.5,8.5);
    \draw[dashed,line width=0.1mm] (54.5,-5.5) --(54.5,8.5);
    \draw[solid,line width=0.4mm] (\a,\a)-- ++\Nu--++\Ed--++\Nu--++\Nd--++\Eu--++\Nd--++\Nu--++\Ed--++\Eu--++\Nd--++\Eu--++\Nd--++\Eu--++\Nd--++\Nu--++\Ed--++\Nu--++\Nd--++\Eu;
    \filldraw (\a,\a) circle (8pt);
    \filldraw (\a,\a)++\Nu circle (8pt);
    \filldraw (\a,\a)++\Nu++\Ed circle (8pt);
    \filldraw (\a,\a)++\Nu++\Ed++\Nu circle (8pt);
    \filldraw (\a,\a)++\Nu++\Ed++\Nu++\Nd circle (8pt);
    \filldraw (\a,\a)++\Nu++\Ed++\Nu++\Nd++\Eu circle (8pt);
    \filldraw (\a,\a)++\Nu++\Ed++\Nu++\Nd++\Eu++\Nd circle (8pt);
    \filldraw (\a,\a)++\Nu++\Ed++\Nu++\Nd++\Eu++\Nd++\Nu circle (8pt);
    \filldraw (\a,\a)++\Nu++\Ed++\Nu++\Nd++\Eu++\Nd++\Nu++\Ed circle (8pt);
    \filldraw (\a,\a)++\Nu++\Ed++\Nu++\Nd++\Eu++\Nd++\Nu++\Ed++\Eu circle (8pt);
    \filldraw (\a,\a)++\Nu++\Ed++\Nu++\Nd++\Eu++\Nd++\Nu++\Ed++\Eu++\Nd circle (8pt);
    \filldraw (\a,\a)++\Nu++\Ed++\Nu++\Nd++\Eu++\Nd++\Nu++\Ed++\Eu++\Nd++\Eu circle (8pt);
    \filldraw (\a,\a)++\Nu++\Ed++\Nu++\Nd++\Eu++\Nd++\Nu++\Ed++\Eu++\Nd++\Eu++\Nd circle (8pt);
    \filldraw (\a,\a)++\Nu++\Ed++\Nu++\Nd++\Eu++\Nd++\Nu++\Ed++\Eu++\Nd++\Eu++\Nd++\Eu circle (8pt);
    \filldraw (\a,\a)++\Nu++\Ed++\Nu++\Nd++\Eu++\Nd++\Nu++\Ed++\Eu++\Nd++\Eu++\Nd++\Eu++\Nd circle (8pt);
    \filldraw (\a,\a)++\Nu++\Ed++\Nu++\Nd++\Eu++\Nd++\Nu++\Ed++\Eu++\Nd++\Eu++\Nd++\Eu++\Nd++\Nu circle (8pt);
    \filldraw (\a,\a)++\Nu++\Ed++\Nu++\Nd++\Eu++\Nd++\Nu++\Ed++\Eu++\Nd++\Eu++\Nd++\Eu++\Nd++\Nu ++\Ed circle (8pt);
    \filldraw (\a,\a)++\Nu++\Ed++\Nu++\Nd++\Eu++\Nd++\Nu++\Ed++\Eu++\Nd++\Eu++\Nd++\Eu++\Nd++\Nu ++\Ed++\Nu circle (8pt);
    \filldraw (\a,\a)++\Nu++\Ed++\Nu++\Nd++\Eu++\Nd++\Nu++\Ed++\Eu++\Nd++\Eu++\Nd++\Eu++\Nd++\Nu ++\Ed++\Nu++\Nd circle (8pt);
    \filldraw (\a,\a)++\Nu++\Ed++\Nu++\Nd++\Eu++\Nd++\Nu++\Ed++\Eu++\Nd++\Eu++\Nd++\Eu++\Nd++\Nu ++\Ed++\Nu++\Nd++\Eu circle (8pt);
\end{tikzpicture}
\end{center}
\caption{A grand zigzag knight's path staying between the lines $y=-2$ and $y=+3$, with size 27 and altitude~0.}
    \label{zigzag staying above a linec}
\end{figure}
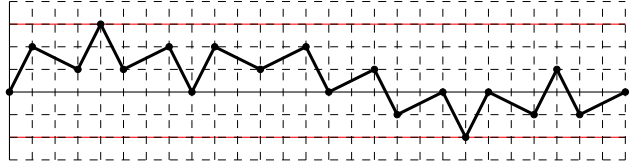

 Here we provide, without going into details, additional results that can be obtained by similar methods. We enumerate grand zigzag knight's paths that stay between the lines $y=-m$ and $y=+M$ for $m,M\geq 0$ with $M\geq m$. In particular we always have $M\geq 1$ (the case $m=M=0$ being trivial). Figure~\ref{zigzag staying above a linec} shows an example of a grand zigzag knight's path staying between the lines $y=-2$ and $y=+3$, with size 27 and altitude $0$.

\begin{rmk}
    When $n=k \Mod{2}$, there is a bijection between those paths ending at $(n,k)$ and pairs of compositions $(X,Y)\in\mathcal{C}_{\frac{n-k}{2},\frac{n+k}{2}}$ such that for all $1\leq j\leq k$,
    $$ -m\leq\sum_{i=1}^j (y_i-x_i)\leq M$$
    (see bijections $\phi$ and $\psi$ in Lemma~\ref{Bijection C Z}).
\end{rmk}

With a similar reasoning as in Section~\ref{bounded grand zigzag}, we obtain the following equations:
\begin{align*}
    F(u)&=u^m-z^3u^{m+M+1}f_M+zu(z+u)(G(u)+u^m),\\
    G(u)&=-\frac{z}{u}f_{-m+1}+\left(\frac{z}{u^2}+\frac{z^2}{u}\right)(F(u)-\mathbb{1}_{[m=0]}).
\end{align*}

\begin{rmk}
    As in Remark \ref{rational solution}, the solutions are rational.
\end{rmk}

Using the kernel method, we find:
$$f_{-m+1}=\scalemath{0.7}{\frac{(1 + sz)\left((\mathbb{1}_{[m=0]} - s^m(1 + sz^2 + s^2z))r^{m + M + 2} + s^{m + M + 1}(r^{m + 1}(1 + rz^2 + r^2z) - \mathbb{1}_{[m=0]}z^2(z + r)(1 + zr))\right)}{s^{m + M}z^2(r + z(2 + sz)) - r^{m + M + 1} }},$$

$$f_M=-\scalemath{0.8}{\frac{r^m + z(z + r)\left(r^{m + 1} - z\left(s^{m-1}(1 + sz)(1 + s^2z + sz^2) - \mathbb{1}_{[m=0]}(r - s)\right)\right)}{z^3(z^2(r + z)(1 + sz)s^{m + M}  - r^{m + M + 1})}}.$$

\begin{theorem}
The generating functions $F(u)$ and $G(u)$ counting the number of grand zigzag knight's paths staying between $y=-m$ and $y=+M$ with respect to the size and the altitude are given by:
\begin{align*}
F(u)&=\scalemath{0.85}{\frac{z^3u^{m + M + 2}f_M + uz^2(u + z)f_{-m+1} + \mathbb{1}_{[m=0]}z^2(u + z)(1 + uz) - u^{m+1}(1+uz^2+u^2z)}{z^3(u-r)(u-s)}},\\
G(u)&=\scalemath{0.85}{\frac{z^3u^{m+M+1}(1 + uz)f_M + uf_{-m+1} + \mathbb{1}_{[m=0]}(1+uz) - u^m(1 + uz)(1 + uz^2 + u^2z)}{uz^2(u-r)(u-s)}}.
\end{align*}
\end{theorem}
We set $H(u):=F(u)+G(u)$.
\begin{corollary}
    The generating function for the number of grand zigzag knight's path staying between $y=0$ and $y=+M$ and ending on the $x$-axis is 
    $H(0)=\frac{f_1}{z^2}.$
\end{corollary}

    When $M=2$, we obtain the first terms, $0\leq n\leq 18$, 
    $$1,\, 0,\, 1,\, 0,\, 2,\, 0,\, 4,\, 0,\, 7,\, 0,\, 14,\, 0,\, 26,\, 0,\, 50,\, 0,\, 95,\, 0,\, 181,$$
    where the even terms are \href{https://oeis.org/A052535}{A052535}  (see Example~\ref{zigzag small tube}).
\begin{rmk}
    We have a convergence to the generating function for all grand  zigzag knight's path staying above the $x$-axis and ending on it $A=\frac{r}{z^3}$ (see \cite[Corollary 1]{BarRa}) as $M\to\infty$.
\end{rmk}

\begin{corollary}
    The generating function for the total number of grand zigzag knight's paths staying between $y=-m$ and $y=+M$ is 
    $$\scalemath{1}{H_{m,M}(z):=H(1)=\frac{(f_{-m+1} + z^2f_M - 1)z +\mathbb{1}_{[m=0]}z(1+z) - 1 - z^2}{z^2 + z - 1}.}$$
\end{corollary}

As a byproduct and using the inclusion–exclusion principle, the generating function for grand zigzag knight's paths staying between lines $y=-m$ and $y=M$, and reaching them, is $H_{m,M}(z)-H_{m-1,M}(z)-H_{m,M-1}(z)+H_{m-1,M-1}(z)$ when $1\leq m\leq M$, and $H_{0,M}(z)-H_{0,M-1}(z)$ when $0\leq M$. Then the generating function grand zigzag knight's paths such that the difference between the maximal and minimal $y$-coordinates of a point of the path is exactly $k$ is $$\scalemath{0.8}{2\left(H_{0,k}(z)-H_{0,k-1}(z) + \sum_{m=1}^{\lfloor k/2\rfloor} \left(H_{m,k-m}(z)-H_{m-1,k-m}(z)-H_{m,k-m-1}(z)+H_{m-1,k-m-1}(z)\right)\right).}$$

Table~\ref{table2} gives the first coefficients for $k=1,2,3$. The sequences do not appear in OEIS.
\begin{table}[ht]
\centering
\rowcolors{2}{}{lightgray}
\begin{tabular}{|c|ccccccccccccccccc|}
    \hline
    \diagbox[]{$k$}{$n$} & 0 & 1 & 2 & 3 & 4 & 5 & 6 & 7 & 8 & 9 & 10 & 11 & 12 & 13 & 14&15&16\\
    \hline
    1 & 0& 0& 2& 0& 2& 0& 2& 0& 2& 0& 2& 0& 2& 0 &2&0&2\\
    2 & 0& 2& 2& 6& 6& 12& 14& 24& 30& 46& 60& 88& 118& 168& 228& 320& 438\\
    3 & 0& 0& 0& 0& 2& 4& 10& 16& 32& 52& 94& 148& 252& 392& 648& 996& 1612\\
    \hline
\end{tabular}
\caption{The number of grand zigzag knight's paths such that the difference of the maximal and minimal heights is at most $k$.}
\label{table2}
\end{table}

\vspace{1em}
\textbf{Acknowledgments.}  The authors would like to thank the anonymous referees for their useful remarks and comments. This work was supported in part by ANR-22-CE48-0002 funded by l’Agence Nationale de la Recherche and by the project ANER ARTICO funded by
Bourgogne-Franche-Comté region (France). José L. Ramírez was partially supported by Universidad Nacional de Colombia, Project No. 57340.

\end{document}